\documentclass[12pt,a4paper]{article}
\usepackage[margin=1in]{geometry}
\usepackage[latin1]{inputenc}
\usepackage{amsmath,bm}
\usepackage{amsthm}
\usepackage{amsfonts}
\usepackage{amssymb}
\usepackage[makeroom]{cancel}
\usepackage{graphicx}
\usepackage{epsfig}
\usepackage{stmaryrd}
\usepackage{fullpage}
\usepackage{color}
\usepackage{booktabs}
\usepackage{colortbl}
\usepackage{hyperref}
\usepackage[noblocks]{authblk}
\usepackage{caption,subcaption}
\usepackage{float}
\usepackage{wrapfig}
\usepackage{caption}
\usepackage{cite}
\usepackage{cleveref}
\usepackage{cases}
\usepackage{enumitem}
\usepackage[titletoc,title]{appendix}
\usepackage{tikz}
\usepackage{ifthen}
\DeclareMathSizes{10}{10}{10}{10}
\title{Capillary Hysteresis and Gravity Segregation in Two Phase Flow Through Porous Media}
\author[1]{K. Mitra\footnote{Corresponding author: email: \href{mailto:koondanibha.mitra@inria.fr}{koondanibha.mitra@inria.fr} }}
\affil[1]{INRIA Paris}
\author[2,3]{C.J. van Duijn}
\affil[2]{Eindhoven University of Technology, Department of  Mechanical Engineering}
\affil[3]{Utrecht University, Department of Earth Sciences}
\date{} 
\numberwithin{equation}{section}
\makeatletter
\DeclareMathSizes{\@xpt}{\@xpt}{6}{5}
\makeatother

\newcounter{assumption}
\newcounter{CaseA}
\newcounter{CaseB}

 \stepcounter{assumption}
 \stepcounter{CaseA}
 \stepcounter{CaseB}

  \newcounter{GSproperties}
 \stepcounter{GSproperties}
 
\newtheorem{remark}{Remark}[section]
\theoremstyle{definition}

\def \a  {\alpha}

\def \d  {\delta}
\def \D  {\Delta}

\def \t  {\tau}

\def \p  {\partial}
\def \R  {\mathbb{R}}

\def \H  {{\cal H}}

\def \S {\mathsf{S}}
\def \P {\mathsf{p}}

\def \Pim {p_c^{(i)}}
\def \Pdr {p_c^{(d)}}

\newcommand{\jump}[1]{ \llbracket #1 \rrbracket}

\begin{document}
\maketitle

\begin{abstract}
We study the gravity driven flow of two fluid phases in a one dimensional homogeneous porous column when history dependence of the pressure difference between the phases (capillary pressure) is taken into account. In the hyperbolic limit, solutions of such systems satisfy the Buckley-Leverett equation with a non-monotone flux function. However,  solutions for the hysteretic case do not converge to the classical solutions in the hyperbolic limit in a wide range of situations. In particular, with Riemann data as initial condition, stationary shocks become possible in addition to classical components such as shocks, rarefaction waves and constant states. We derive an admissibility criterion for the stationary shocks and outline all admissible shocks. Depending  on the capillary pressure functions, flux function and the Riemann data, two cases are identified a priori for which the solution consists of a stationary shock. In the first case, the shock remains at the point where the initial condition is discontinuous. In the second case, the solution is frozen in time in at least one semi-infinite half.  The predictions are verified using numerical results.  
\end{abstract}

\section{Introduction}\label{sec:Intro}
\begin{figure}[h!]
\begin{subfigure}{.5\textwidth}
\includegraphics[scale=.4]{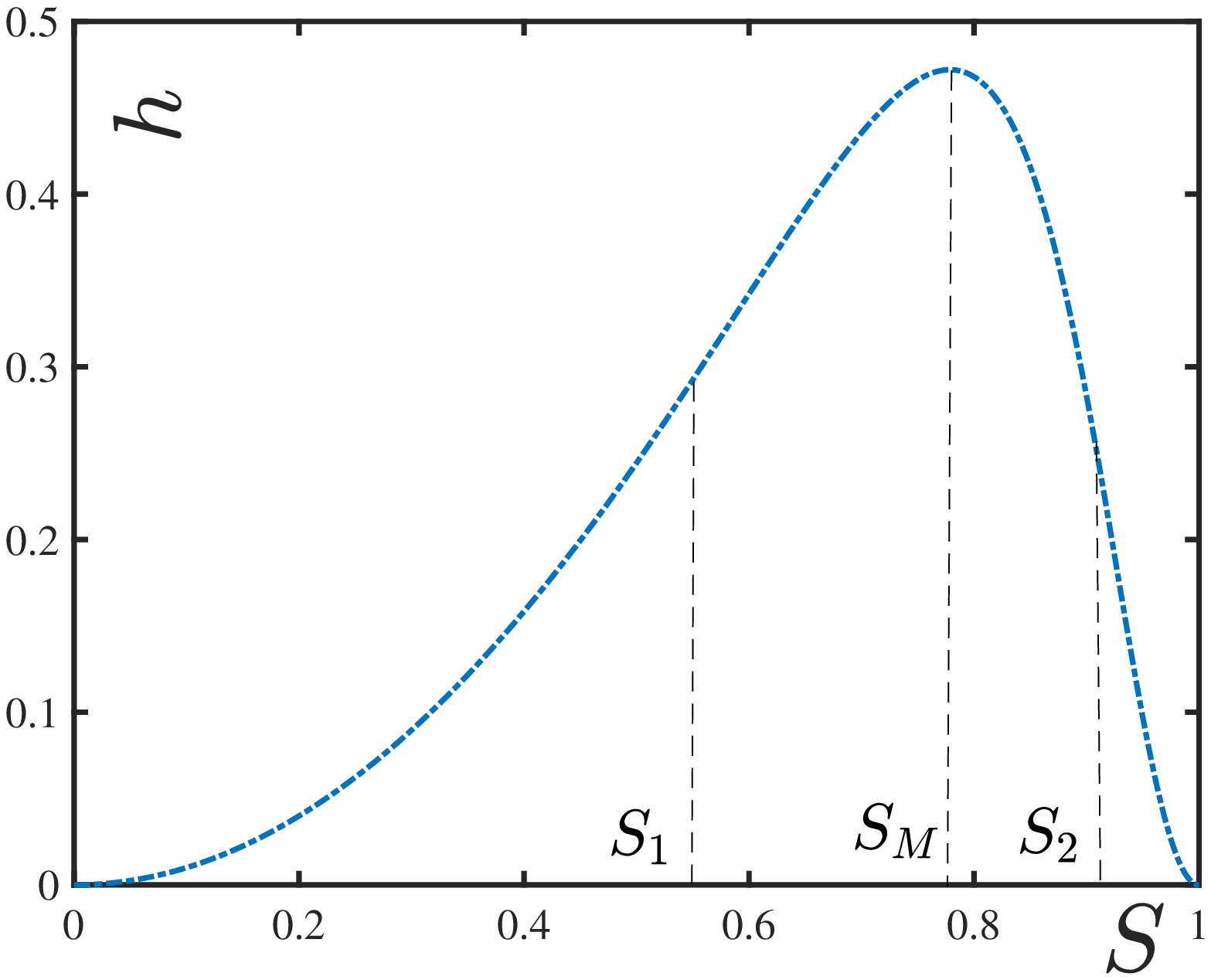}
\end{subfigure}
\begin{subfigure}{.5\textwidth}
\includegraphics[scale=.4]{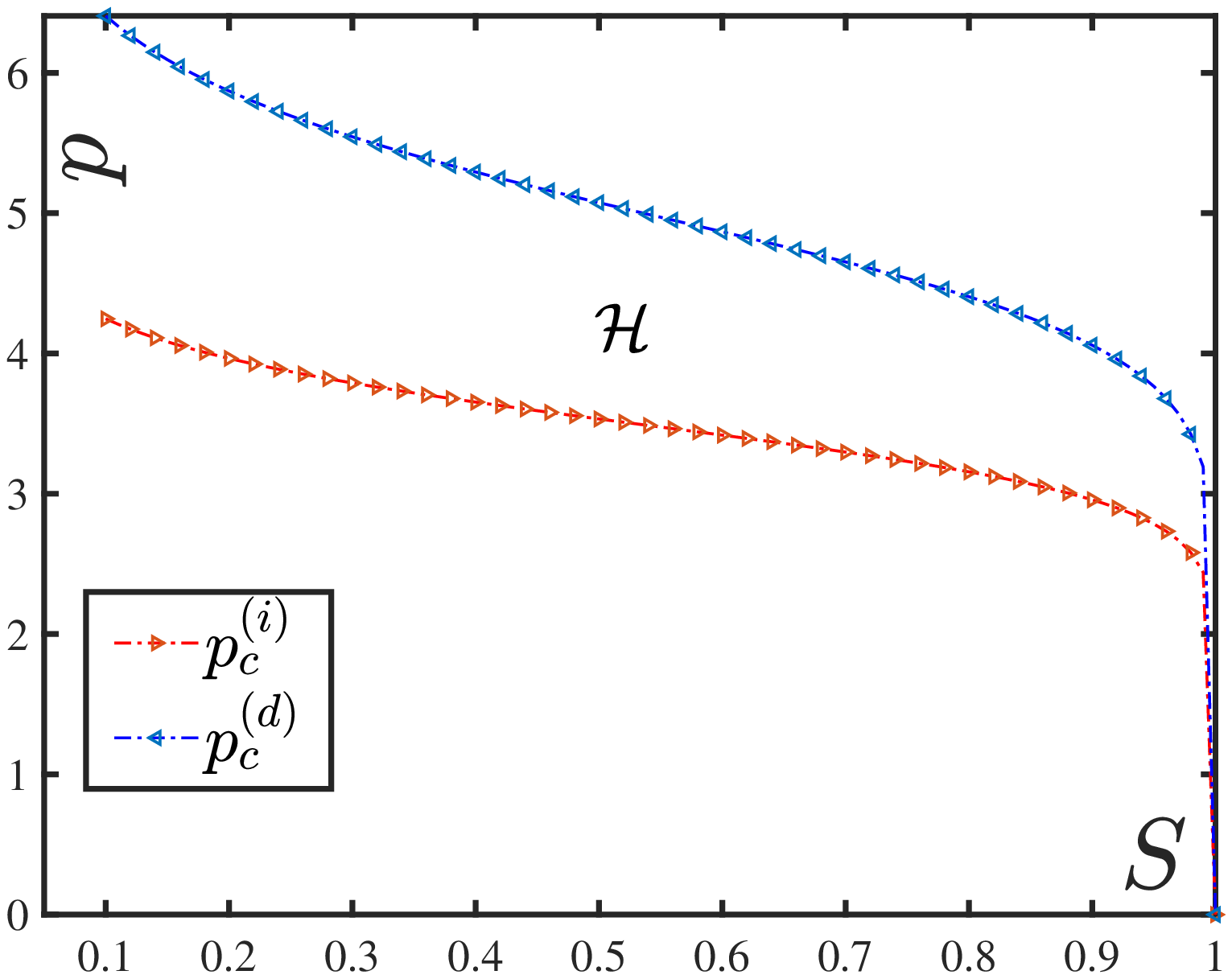}
\end{subfigure}
\caption{(left) The $h$-$S$ curve for water (wetting phase, viscosity $5.23\times 10^{-4} \mathrm{[Pa \cdot s]}$)  and methane (viscosity $1.202\times 10^{-5}  \mathrm{[Pa \cdot s]}$). In reference to \Cref{sec:MathIntro}, the Brooks-Corey relationship has been used for the relative permeabilities, the mobility ratio is $M=43.52$ and the characteristic saturations  $S_M$, $S_1$, $S_2$ are marked. (right) The imbibition and drainage capillary pressure curves and the hysteretic region $\H$. The van Genuchten model has been used for $p_c^{(i)}$ and $\Pdr$. Details are given in the numerical section. }
\label{fig:Pch}
\end{figure}

In this paper we investigate gravity driven flow of two fluid phases through a homogeneous one-dimensional porous column. 
We are concerned with the special case in which the length of the column is large and no injection of fluid is present (i.e. the total flow is zero). The corresponding mathematical model uses the well-known 
 Buckley-Leverett equation, see \cite{helmig1997multiphase}. In dimensionless form it reads
\begin{equation}
\p_t S + \p_x h(S)=0.\label{eq:HBL}
\end{equation}
Here $S\in [0,1]$ is the saturation of the wetting phase and $h:[0,1]\to [0,\infty)$ the fractional flow function of which a typical sketch is shown in \Cref{fig:Pch} (left). The space coordinate is  $x$ and $t$ denotes time. We solve \eqref{eq:HBL} for $t>0$ and $x\in \R$, where we prescribe at $t=0$ the Riemann condition
\begin{equation}\label{eq:iniS}
S(x,0)=\begin{cases}
S_T &\text{ for }  x<0,\\ 
S_B &\text{ for }  x>0, 
\end{cases} \text{ where } 0<S_B<S_T<1 \text{ are constants}.
\end{equation}

Solutions of \eqref{eq:HBL}--\eqref{eq:iniS} are generally non-unique \cite[Chapter 1]{lefloch2002hyperbolic}. To find the physically relevant solution, we replace \eqref{eq:HBL} by the capillary Buckley-Leverett equation
\begin{equation}
\p_t S + \p_x [ h(S)(1+ \d\, \p_x p)]=0.\label{eq:BL}
\end{equation}
and study the limit $\d\to 0$. In \eqref{eq:BL}, $p$ is the pressure difference between the fluid phases and  $\d$ is the dimensionless capillary number which scales inversely with the length of the domain (reference length). Hence for large domains, \eqref{eq:HBL} is approximated by \eqref{eq:BL}. A detailed derivation is given in \Cref{sec:MathIntro}.

To solve \eqref{eq:BL} a relation between $S$ and $p$ is assumed. In the standard equilibrium approach, one uses
\begin{align}\label{eq:StandardPcS}
 p=p_c(S),
\end{align}
where $p_c: (0,1)\to [0,\infty)$ is a capillary pressure function \cite{helmig1997multiphase}. 

For $\d>0$, let $(S_\d,p_\d)$ denote the solution of \eqref{eq:iniS}--\eqref{eq:StandardPcS}. The hyperbolic limit as $\d\to 0$ yields the well known Buckley-Leverett solution of \eqref{eq:HBL}--\eqref{eq:iniS}, comprising of constant states separated by shocks and  rarefaction waves \cite{oleinik1957discontinuous,liu2000hyperbolic,lefloch2002hyperbolic}. In particular, shocks can be seen as the limit of smooth, monotone travelling waves which become steeper as $\d\to 0$ \cite{lefloch2002hyperbolic,liu2000hyperbolic,van1990large}. It is also known that the hyperbolic limit solution is independent of the actual shape of the capillary pressure $p_c$ as long as the approximating equation \eqref{eq:BL} is of convection-diffusion type \cite[Chapter 3]{lefloch2002hyperbolic}.  

However, there are other more realistic capillary pressure expressions for which the limiting $(\d\to 0)$ solution does inherit some properties of the vanishing capillarity. This was studied in detail in \cite{van2007new,van2013travelling,spayd2011buckley,mitra2019fronts} for the case where \eqref{eq:StandardPcS} is replaced by the non-equilibrium expression
$$
p=p_c(S)-\t \p_t S,
$$
where $\t>0$ is a relaxation parameter called dynamic capillarity coefficient which attributes to saturation overshoot \cite{VANDUIJN2018232,mitra2018wetting,MyThesis}.

In this paper we investigate the effect of hysteresis in the capillary pressure. Since $p_c(\cdot)$ is a single valued function of saturation only, it does not contain any information about the history or directionality of the process. In particular, it does not distinguish between imbibition and drainage. No hysteretic effects are present in  \eqref{eq:StandardPcS}. However, hysteresis is known to occur in multi-phase porous media flow. This was first observed by Haines in 1930 and has been verified since by numerous experiments, some notable examples being \cite{Morrow_Harris,zhuang2017advanced}. An overview of different hysteresis models from the mathematical, modelling and physical perspective can be found in \cite{schweizer2017hysteresis,Kmitra2017,miller2019nonhysteretic}.  

In our approach we replace \eqref{eq:StandardPcS} by the following hysteresis description. Let
$$\Pim,\, \Pdr:(0,1)\to (0,\infty), \quad \Pim(S)< \Pdr(S) \text{ for } 0<S<1,$$
denote the imbibition and drainage capillary pressure functions, typical examples being shown in \Cref{fig:Pch} (right),
 and let
$$
\H:= \{(S,p): 0\leq S\leq 1,\; \Pim(S)\leq p\leq \Pdr(S)\}
$$  
be the hysteresis region in the $(S,p)$-plane, as indicated in \Cref{fig:Pch} (right). We restrict ourselves to the well-known play-type hysteresis model, first proposed in \cite{Beliaev2001}, which relates $S$ and $p$ via the relation
\begin{align}\label{eq:play-type}
 p\in \H, \text{ and } p=\begin{cases}
\Pim(S) &\text{ when } \p_t S> 0,\\
\in [\Pim(S),\Pdr(S)] &\text{ when } \p_t S= 0,\\
  \Pdr(S) &\text{ when } \p_t S< 0.
 \end{cases}
 \end{align} 
The play-type model assumes that switching from drainage to imbibition, or vice versa, occurs only along vertical scanning curves. For the rest of the study the state when $\p_t S=0$ and consequently $p\in [\Pim(S),\Pdr(S)]$ will be referred to as the undetermined state of the system, whereas $\p_t S>0$ (consequently $p=\Pim(S)$) and $\p_t S<0$ ($p=\Pdr(S)$) will refer to the imbibition and drainage states, respectively.
Other models, such as the Lenhard-Parker model \cite{parker1987parametric} and the extended play-type hysteresis model \cite{Kmitra2017,KMitraEUX_2020} assume a more complex relation between $S$ and $p$ when $\Pim(S)<p<\Pdr(S)$. We comment on the consequences of such models in \Cref{rem:Generality}.

\begin{figure}[h!]
\begin{subfigure}{.5\textwidth}
\includegraphics[scale=.4, trim={0 0 0 -1cm},clip]{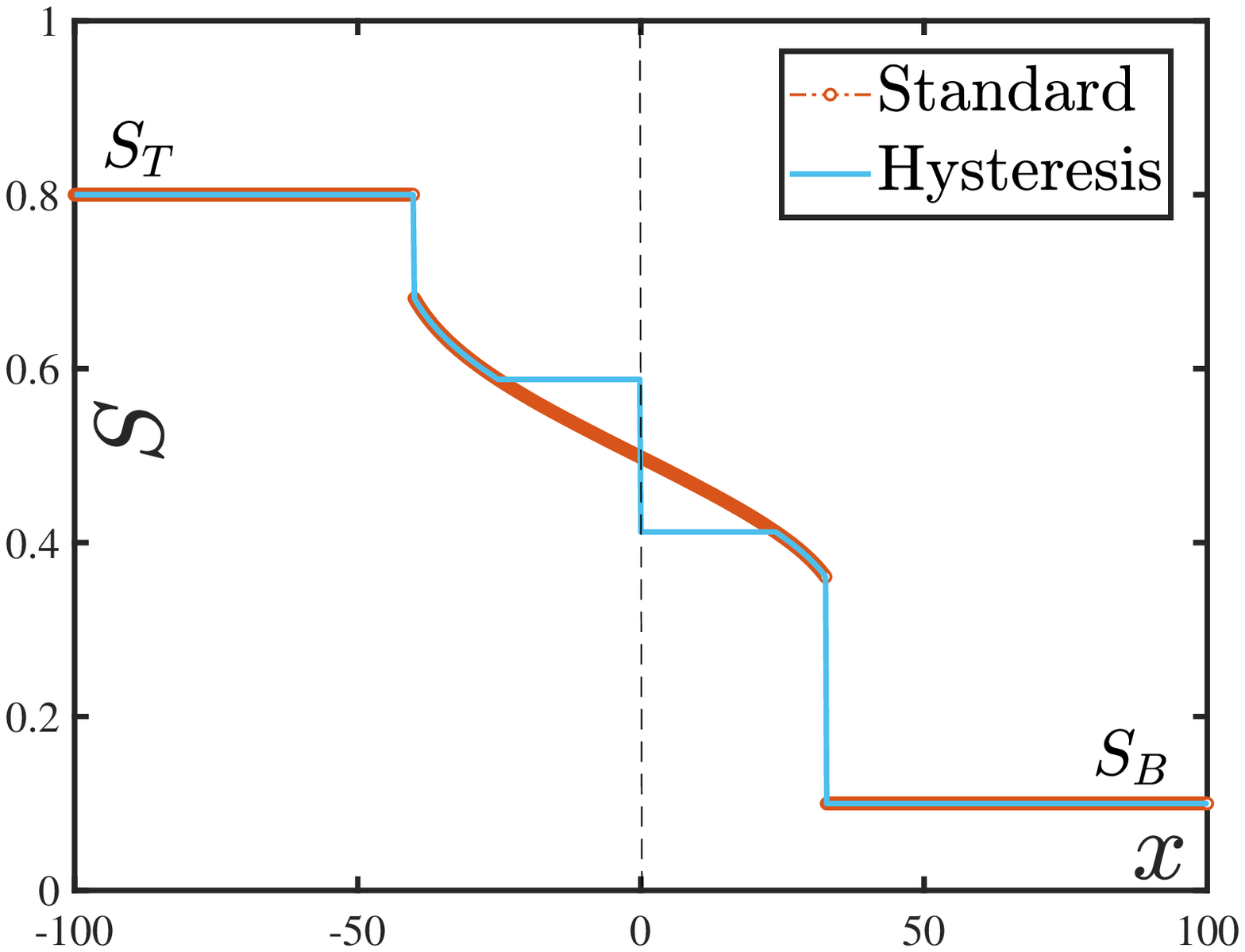}
\end{subfigure}
\begin{subfigure}{.5\textwidth}
\includegraphics[scale=.4, trim={0 0 0 -1cm},clip]{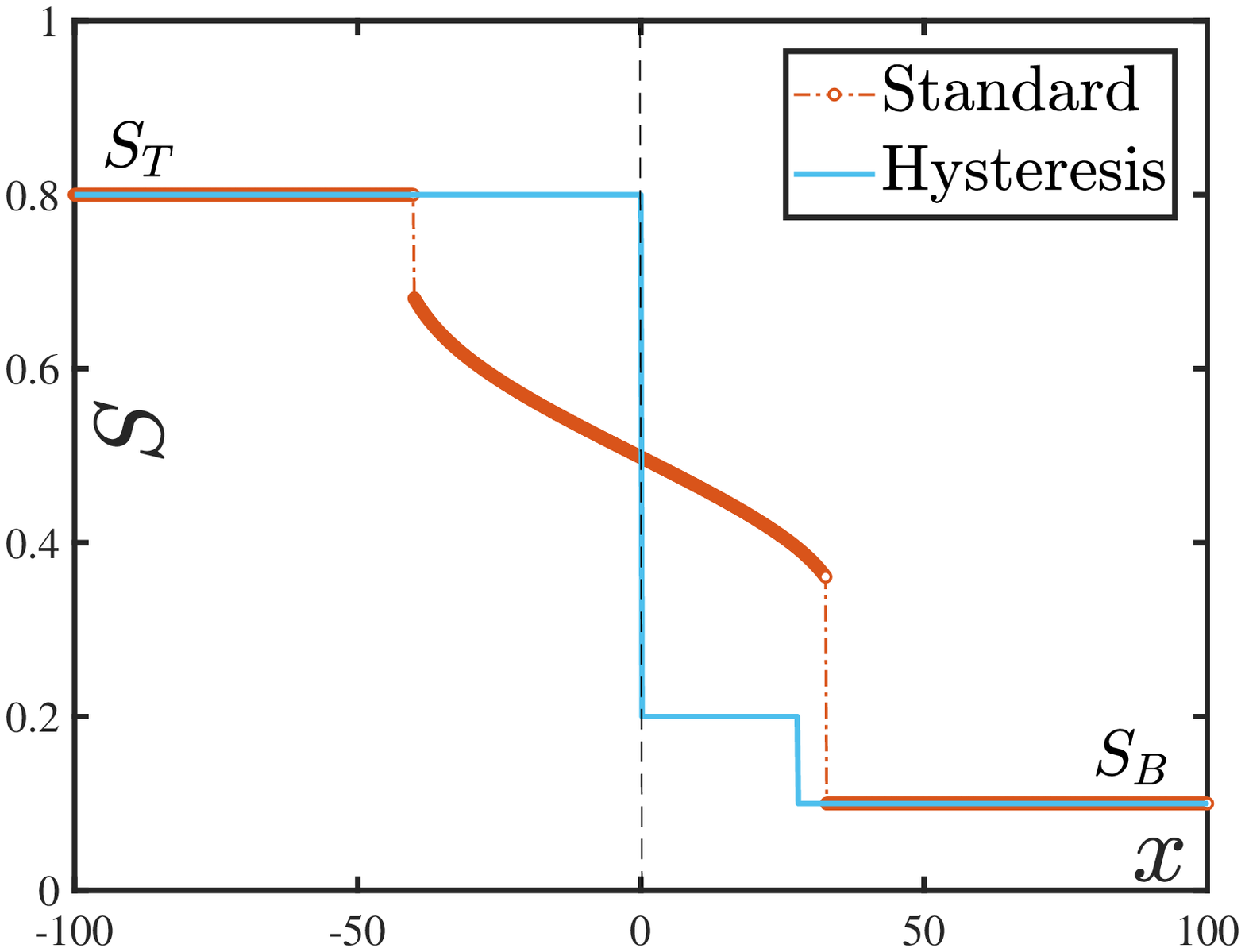}
\end{subfigure}
\caption{Vanishing capillarity solutions of \eqref{eq:HBL} with initial condition \eqref{eq:iniS}. Here, $S_B=0.1$, $S_T=0.8$ and $t=100$. The red marked profiles show the classical solutions (for the standard model \eqref{eq:StandardPcS}), whereas, the blue solid profiles show the result for the play-type model \eqref{eq:play-type}. For the (left) plot, $\Pim$ and $\Pdr$ curves are close to each other, whereas in the (right) plot the curves are as shown in \Cref{fig:Pch} (right). Exact details are given later. }\label{fig:StandardVsHys}
\end{figure}

The purpose of this paper is to construct solutions of the Riemann problem \eqref{eq:HBL}--\eqref{eq:iniS} that arise as the hyperbolic limit ($\d\to 0$) of \eqref{eq:iniS}--\eqref{eq:BL} and \eqref{eq:play-type}. We demonstrate that the occurrence of hysteresis in the vanishing capillary term, i.e. using \eqref{eq:play-type} instead of \eqref{eq:StandardPcS}, gives solutions that are significantly different specially when $S_B$ is close to $0$ and $S_T$ is close to $1$. In particular, a stationary discontinuity occurs at the location where the initial condition is discontinuous. This is illustrated in \Cref{fig:StandardVsHys}. 
The magnitude of the jump depends on the difference $\Pdr(S)-\Pim(S)$, as will be discussed later.
\noindent
This behaviour was mentioned briefly in Section 3.5 of \cite{mitra2019fronts}. A stationary discontinuity in saturation for the classical problem (with no hysteresis) was also predicted in \cite{andreianov2013vanishing} at the interface between two semi-infinite halves having different $h$ and $p_c$ functions. However, the directionality imposed by hysteresis demands an extension of their results.  The saturation discontinuity has major practical importance as the saturation distribution can change considerably if hysteresis is present, as is evident from \Cref{fig:StandardVsHys}.

\begin{remark}[Vanishing capillarity method]
In the mathematical literature, the method of finding solutions of equations such as \eqref{eq:HBL} by passing to the limit $\d\to 0$ in \eqref{eq:BL} is called the vanishing viscosity method, and was studied in classical references such as \cite{lax1973hyperbolic,oleinik1957discontinuous}.  In these papers, the notion of entropy was introduced, calling the vanishing viscosity solutions `entropy' solutions. In the context of our application, we call the approach vanishing capillarity method, the approximate solutions $(S_\d,p_\d)$ for $\d>0$ the capillarity solutions, and the  $\d\to 0$ limit solution the vanishing capillarity solution.
\end{remark}

Hysteresis models have been analysed exhaustively, particularly in one spatial dimension \cite{MyThesis}. Existence  and uniqueness  results for the regularised play-type model are given in \cite{lamacz2011well,SCHWEIZER20125594,koch2013,KMitraEUX_2020} and \cite{Cao2015688} respectively. A horizontal redistribution study using similarity solutions was performed in \cite{Kmitra2017}. Travelling wave analysis for hysteresis was conducted for flow of water through soil in \cite{VANDUIJN2018232,mitra2018wetting,el2018traveling}. For the two-phase case, travelling waves were studied in \cite{mitra2019fronts} for monotone flux functions and vanishing capillarity solutions were derived that differ significantly from the classical ones. Non-classical solutions for non-monotone flux functions like $h$ are investigated in \cite{shearer2015traveling}, although hysteresis is not included. In \cite{krejci1996hysteresis}, examples of how hysteresis influences hyperbolic solutions in mechanics are found.  The role of relative permeability hysteresis, which is not addressed in this study, in determining the hyperbolic solutions is examined in \cite{plohr2001modeling,schaerer2006permeability,bedrikovetsky1996mathematical,abreu2012computational}. 
However, the relation between $p$ and $S$ is assumed to be a linear one in these articles. In the current study, we investigate the vanishing capillarity solutions for the capillary hysteresis models in the non-monotone flux case and show that non-classical behaviour such as stationary shocks may occur. Stationary discontinuities have been studied in \cite{andreianov2013vanishing} for heterogeneous media  without hysteresis and for redistribution problems (no gravity) in \cite{Kmitra2017,Philip_1991} with hysteresis. The occurrence of a stationary discontinuity due to hysteresis in gravity driven flows is to our knowledge a novel observation.

We structure the paper as follows: In \Cref{sec:MathIntro} the assumptions are stated and the model is derived. Using travelling wave analysis,  all admissible shocks, including stationary shocks, are derived in \Cref{sec:Admissible} for the standard and the hysteresis model. Then, in \Cref{sec:EntropySol}, the admissible shocks are used to construct the vanishing capillarity solutions. Two cases are identified when the solutions for hysteresis deviate from the classical ones. This is determined a priori from  $h$, $\Pim$, $\Pdr$, $S_B$ and $S_T$. In the first case, the solution has a stationary discontinuity while the rest of the solution retains the structure of the classical solutions (\Cref{sec:Case1}). The second case has no classical counterparts and the  solution is frozen in one of the semi-infinite halves (\Cref{sec:Case2}). Finally, in \Cref{sec:Num}, we solve \eqref{eq:iniS}--\eqref{eq:BL} and \eqref{eq:play-type} numerically  for small $\d>0$ and show that the solution closely resembles our predictions.

\section{Problem Formulation}\label{sec:MathIntro}
We set-up the two-phase flow problem in a one-dimensional homogeneous porous domain. The phases are assumed to be incompressible and immiscible. There is no injection at the boundaries, making the flow purely counter-current and gravity driven, having zero total flux of the combined wetting and non-wetting phases. Following \cite{helmig1997multiphase}, we consider for each phase the mass balance equation and the corresponding Darcy Law. This yields
\begin{subequations}\label{eq:TwoPhase}
\begin{align}
\phi\,\p_t S=&\p_x\left [\frac{K }{\mu_w}k_w(S)(\p_x p_w -\rho_w g)\right ] \qquad &\text{(wetting phase)},\label{eq:Weq}\\
\phi\, \p_t (1 - S)=&\p_x\left [\frac{K }{\mu_n}k_n(1-S)(\p_x p_n -\rho_n g)\right ]\qquad &\text{(non-wetting phase)}\label{eq:Neq}.
\end{align}
\end{subequations}
For each phase $\a=w,n$ ($w$ and $n$ representing the wetting and the non-wetting phases respectively), $p_\a$ denotes the pressure, $k_\a$ the relative permeability, $\mu_\a$ the viscosity and $\rho_\a$ the density. The porosity $\phi$ and absolute permeability $K$ are properties of the medium and are constant due to the assumption of homogeneity. Finally, gravity points in the direction of positive $x$ and $g$ is the gravitational constant.
 For the remainder of the study we assume that the wetting phase is denser than the non-wetting phase, i.e.,
\begin{equation*}
\rho_n<\rho_w.
\end{equation*}
 Adding the equations in \eqref{eq:TwoPhase} gives
$$
\p_x \left[\tfrac{K k_w}{\mu_w}(\p_x p_w -\rho_w g) + \tfrac{K k_n}{\mu_n}(\p_x p_n -\rho_n g)\right ]=0,
$$
The term inside the brackets $[\,]$ is the total flux of the  combined phases. Since no-injection takes place in the column
$$
\tfrac{K k_w}{\mu_w}(\p_x p_w -\rho_w g) + \tfrac{K k_n}{\mu_n}(\p_x p_n -\rho_n g)=0.
$$
Defining the capillary pressure and  mobility ratio as
$$
p:=p_n-p_w \quad \text{ and } \quad M:=\mu_w\slash \mu_n
$$ respectively, one obtains through rearrangement
\begin{equation}
\frac{K k_w}{\mu_w}(\p_x p_w -\rho_w g)= -\dfrac{K}{\mu_w} \dfrac{k_n k_w}{k_n + M^{-1} k_w }\p_x p- \dfrac{K g (\rho_w-\rho_n)}{\mu_w} \dfrac{k_n k_w}{k_n + M^{-1} k_w }.\label{eq:WaterFlux}
\end{equation}
To non-dimensionalise \eqref{eq:TwoPhase}, we introduce a characteristic pressure $p_{\!_{\mathrm{ref}}}$ (taken from the capillary pressure curves), a characteristic length $H$ (length of column or typical observation distance) and a characteristic time $t_{\!_{\mathrm{ref}}}$. We further define the fractional flow function
\begin{equation}
h:=\dfrac{k_n k_w}{k_n + M^{-1} k_w },
\end{equation}
and introduce the dimensionless capillary number
\begin{equation}
\d:=\frac{p_{\!_{\mathrm{ref}}}}{(\rho_w-\rho_n) g H}>0.
\end{equation}
Inserting \eqref{eq:WaterFlux} into \eqref{eq:Weq}, choosing $t_{\!_{\mathrm{ref}}}=\frac{\mu_w \phi H}{K g (\rho_w-\rho_n)}$, and redefining the dimensional quantities as their dimensionless versions
$$
p\mapsto \tfrac{p}{p_{\!_{\mathrm{ref}}}}, \quad t\mapsto \tfrac{t}{t_{\!_{\mathrm{ref}}}}, \quad x\mapsto \tfrac{x}{H},
$$
one obtains \eqref{eq:BL}. 
The Brooks-Corey model is commonly used for the relative permeability functions:
\begin{equation}
k_\a(S)=S^2 \quad \text{ for } \a=n,w. \label{eq:BrooksCorey}
\end{equation}
This gives the shape of $h$ as shown in \Cref{fig:Pch} and properties outlined in \ref{ass:h}. As for $p$, either \eqref{eq:StandardPcS} or \eqref{eq:play-type} is used.

For the capillary pressures and the fractional flow function, following set of properties is assumed, see \Cref{fig:Pch}. They are consistent with experimental observations \cite{helmig1997multiphase,Bear1979}:
\begin{enumerate}[label=(P\theGSproperties)]
  \item The capillary pressures $p_c^{\left(j \right)}: \left(0,1 \right] \rightarrow \left(-\infty,\infty \right)$, for $j=i,d$,  are continuously differentiable and strictly decreasing in $(0,1)$;  
  $p_c^{\left( i \right)}\left(S \right) < p_c^{\left( d \right)}\left(S \right) \text{ for } S \in \left(0,1 \right)$ and $\Pim(1)=\Pdr(1)= 0$ (no entry pressure).  \label{ass:Pc}
\stepcounter{GSproperties}
 \item The fractional flow function $h:[0,1]\to [0,\infty)$ is smooth with $h(0)=h(1)=0$. There exists $S_M\in (0,1)$ such that 
$$
h'(S)>0 \text{ for } 0<S<S_M, \text{ and } h'(S)<0 \text{ for } S_M<S<1.
$$
Moreover, there exists inflection points $S_1,\,S_2\in (0,1)$ with $0<S_1<S_M<S_2<1$ such that 
$$
h''(S)>0 \text{ for } \{0<S<S_1\}\cup \{S_2<S<1\}, \text{ and }  h''(S)<0 \text{ for } S_1<S<S_2.
$$
 \stepcounter{GSproperties}\label{ass:h}
\end{enumerate}
\vspace{-1em} 
The assumption \ref{ass:Pc} is consistent with the van Genuchten model for capillary pressures.
For $h$, the Brooks-Corey model \eqref{eq:BrooksCorey} yields
\begin{equation}\label{eq:hcurve}
h(S)=\frac{S^2(1-S)^2}{M ^{-1} S^2 + (1-S)^2} \quad \text{ with } \quad S_M=\frac{1}{(1+ M^{-1/3})}.
\end{equation}

\section{Admissible shocks}\label{sec:Admissible}
In this section we consider shock solutions of equation \eqref{eq:HBL} that originate from smooth solutions of \eqref{eq:BL}. A shock is characterized by a constant left state $S_l$, a constant right state $S_r$ and a constant speed $c$. They are denoted by $\{S_l,\,S_r,\,c\}$ or
\begin{align}\label{eq:DefShocks}
S(x,t)=\begin{cases}
S_l &\text{ for } x<ct,\\
S_r &\text{ for } x>ct.
\end{cases}
\end{align}
We call  $\{S_l,\,S_r,\,c\}$ an admissible (i.e. vanishing capillarity) shock if it can be approximated, as $\d\to 0$, by smooth solutions $S_\d$ of \eqref{eq:BL} \cite{liu2000hyperbolic,lefloch2002hyperbolic}. To investigate this, we consider a special class of solutions of \eqref{eq:BL} in the form of travelling waves:
\begin{align}
S_\d(x,t)=\S(\eta), \quad p_\d(x,t)=\P(\eta), \text{ with } \eta=\frac{x-ct}{\d}.\label{eq:TravellingWaves}
\end{align}
Here $\S:\R\to [0,1]$ is the saturation profile, $\P:\R\to [0,\infty]$ the capillary pressure profile and $c\in \R$ the wave-speed. We consider profiles that satisfy for the saturation 
\begin{align}\label{eq:TWBCS}
\lim\limits_{\eta\to -\infty} \S(\eta)=S_l, \quad \lim\limits_{\eta\to +\infty} \S(\eta)=S_r, 
\end{align}
and for the pressure 
\begin{equation}\label{eq:TWBCP}
\P\in L^\infty(\R)\quad \text{(bounded pressure)}.
\end{equation}
Clearly, if a smooth profile \eqref{eq:TravellingWaves} exists and satisfies \eqref{eq:TWBCS}, then
$$
S(x,t)=\lim\limits_{\d\to 0} S_\d(x,t)= \lim\limits_{\d\to 0} \S\left (\frac{x-ct}{\d}\right )=\begin{cases}
S_l &\text{ for } x<ct,\\
S_r &\text{ for } x>ct.
\end{cases}
$$
Substituting \eqref{eq:TravellingWaves} into \eqref{eq:BL} gives
\begin{align}
-c\S' + [h(\S)(1+\P')]=0 \text{ in } \R,
\end{align}
where primes denote differentiation. Integrating this expression yields 
\begin{align}
-c\S + h(\S)(1+ \P')=A=\text{ constant in } \R.\label{eq:IntegratedBL}
\end{align}
If $\S$ satisfies \eqref{eq:TWBCS}, then \eqref{eq:IntegratedBL} implies that $\P'$ has a limit for $\eta\to \pm\infty$. The boundedness of pressure \eqref{eq:TWBCP} then forces
\begin{align}
\lim\limits_{\eta\to -\infty}\P'(\eta)=\lim\limits_{\eta\to +\infty} \P'(\eta)=0. 
\end{align}
Using this  in \eqref{eq:IntegratedBL}  for $\eta=\pm\infty$ results in
\begin{subequations}\label{eq:2equations2c}
\begin{align}
&-c S_l + h(S_l)=A,\\
&-c S_r + h(S_r)=A,
\end{align}
\end{subequations}
giving
\begin{align}\label{eq:RankineHugoniot}
c=c(S_l,S_r):=          \frac{h(S_l)-h(S_r)}{S_l-S_r}=:\frac{\jump{h}}{\jump{\S}} \quad \text{  (Rankine-Hugoniot)}.
\end{align}
Substituting $c$ and $A$ in \eqref{eq:IntegratedBL} yields
\begin{align}\label{eq:MainOde}
h(\S)\, \P'=(\S-S_l) [\,c(S_l,S_r)-c(S_l,\S)],\quad  \text{ in } \R.
\end{align}

\subsection{Equilibrium case: $p$ given by \eqref{eq:StandardPcS}}\label{sec:ShocksPcS}
Then $\P=p_c(\S)$, and we write \eqref{eq:MainOde} as
\begin{align}
D(\S)\,\S'=F(\S) \text{ in } \R,\label{eq:ODES}
\end{align}
where 
$$
D(\S)=-h(\S)\frac{d p_c}{dS}(\S)>0, \text{ and } F(\S)=(S_l-\S) [\, c(S_l,S_r)-c(S_l,\S)].
$$
Clearly $F(S_l)=F(S_r)=0$.

\begin{figure}
\centering
\includegraphics[scale=.5]{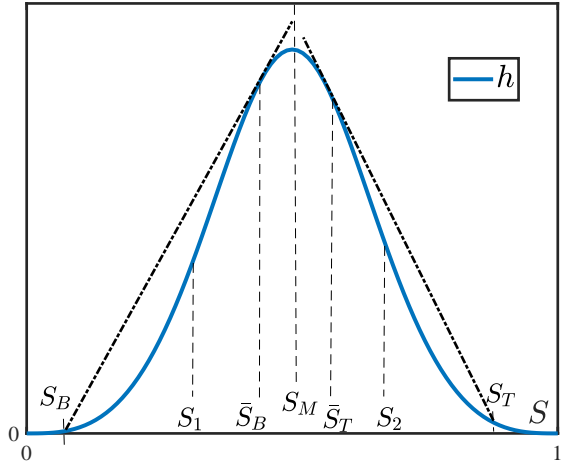}
\caption{Graphical interpretation of conditions \eqref{eq:LaxCondition}, applied to $S_B<S_1$ as the right state and $S_T>S_2$ as the left state. The saturations $\bar{S}_B$ and $\bar{S}_T$ are shown in the figure.}\label{fig:barSbSt}
\end{figure}

If $S_l>S_r$, equation \eqref{eq:ODES} has a solution satisfying \eqref{eq:TWBCS} 
if $F(s)<0$ for all $S_r<s<S_l$. Hence, the shock $\{S_l,\,S_r,\,c\}$ is admissible if $c$ is given by \eqref{eq:RankineHugoniot} and
\begin{subequations}\label{eq:LaxCondition}
\begin{align}
c(S_l,S_r)<c(S_l,S) \text{ for all } S_r<S<S_l.
\end{align}
If $S_l<S_r$, the shock $\{S_l,\,S_r,\,c\}$ is admissible if \eqref{eq:RankineHugoniot} is satisfied, and
\begin{align}
c(S_l,S_r)>c(S_l,S) \text{ for all } S_l<S<S_r.
\end{align}
\end{subequations}
Conditions \eqref{eq:LaxCondition} are classical (Oleinik admissibility conditions), see \cite{oleinik1957discontinuous,liu2000hyperbolic,lefloch2002hyperbolic}. To see which possible admissible shocks connect to the left or right states of the Riemann problem \eqref{eq:HBL}--\eqref{eq:iniS}  we introduce two saturations.
For $S_B<S_1$ (first inflection point of $h$), let $\bar{S}_B\in (S_B,1)$ denote the unique point where the line through $(S_B,h(S_B))$ is tangent to the graph of $h$. If $S_B\geq S_1$, we set $\bar{S}_B=S_B$. A similar definition gives $\bar{S}_T\in (0,S_T)$, see \Cref{fig:barSbSt}. Thus, the points $\bar{S}_B$ and $\bar{S}_T$ satisfy,
\begin{equation}
h'(\bar{S}_B)=\frac{h(\bar{S}_B)-h(S_B)}{\bar{S}_B -S_B}=c(S_B,\bar{S}_B), \quad h'(\bar{S}_T)=\frac{h(\bar{S}_T)-h(S_T)}{\bar{S}_T -S_T}=c(\bar{S}_T,S_T).\label{eq:barS}
\end{equation}
Applying \eqref{eq:LaxCondition} to Riemann problem \eqref{eq:HBL}--\eqref{eq:iniS}, we note that if $S_B<S_1$ , then $S_B$ can serve as the right state $S_r=S_B$ of an admissible shock for left states between $S_B<S_l<\bar{S}_B$, see \Cref{fig:barSbSt}. Similar for $S_l=S_T>S_2$ and  $\bar{S}_T<S_r<S_T$.

\subsection{Hysteretic case: $p$ given by \eqref{eq:play-type}}\label{sec:ShocksPlayType}
If $c(S_l,S_r)>0$ or $c(S_l,S_r)<0$, then any admissible shock in terms of \eqref{eq:StandardPcS} is also an admissible shock in terms of \eqref{eq:play-type}: To see this, fix $S_l>S_r$ and let \eqref{eq:LaxCondition} hold. Define the relation between $\P(\eta)$ and $\S(\eta)$ as
$$
\P(\eta)=\begin{cases}
\Pim(\S(\eta)) &\text{ if } c(S_l,S_r)>0,\\
\Pdr(\S(\eta)) &\text{ if } c(S_l,S_r)<0.
\end{cases}
$$
Then, separately for $c(S_l,S_r)>0$ and $c(S_l,S_r)<0$,  $\S$ and $\P$ are related through the classical relation \eqref{eq:StandardPcS} where $p_c$ is replaced by $\Pim$ and $\Pdr$ respectively. Consequently, the existence of $\S(\eta)$ and $\P(\eta)$ solving \eqref{eq:MainOde} with boundary conditions \eqref{eq:TWBCS}--\eqref{eq:TWBCP} follows from \Cref{sec:ShocksPcS}. 
Observe that $\S'(\eta)<0$ for all $\eta\in \R$ as a consequence of \eqref{eq:ODES} and \eqref{eq:LaxCondition}. Hence, recalling from \eqref{eq:TravellingWaves} that $S_\d(x,t)=\S(\eta)$ and $p_\d(x,t)=\P(\eta)$, we have
\begin{align*}
\p_t S_\d(x,t)=-\frac{c(S_l,S_r)}{\d}\S'(\eta)\begin{cases}
>0 &\text{ when } c(S_l,S_r)>0 \text{ (imbibition wave)},\\
<0 &\text{ when } c(S_l,S_r)<0 \text{ (drainage  wave)}.
\end{cases}
\end{align*}
It is now straightforward to verify that $(S_\d,p_\d)$ satisfies the hysteresis relation \eqref{eq:play-type} for both $c(S_l,S_r)>0$ and $c(S_l,S_r)<0$.  Hence, shocks that are admissible in terms of the equilibrium capillary pressure \eqref{eq:StandardPcS} are admissible in terms of the hysteretic capillary pressure \eqref{eq:play-type} as well, provided $c\not=0$. Note that the profiles $(\S,\P)$ depend on the functional form of $p_c$, but not the resulting shock. 

Observe that for $c(S_l,S_r)>0$, the entire approximating wave $(S_\d,p_\d)$ is in imbibition state since $p_\d(x,t)=\Pim(S_\d(x,t))$ and $\p_t S_\d> 0$. Hence, we also refer to the resulting shock as being in \textbf{imbibition state}. Similarly, for $c(S_l,S_r)<0$, the shock is in \textbf{drainage state} since $p_\d(x,t)=\Pdr(S_\d(x,t))$ and $\p_t S_\d< 0$. 
The effects of hysteresis are only observed, if the hysteretic states ahead and behind the shock are different. This can only happen for stationary shocks having zero wave speed, i.e. $c(S_l,S_r)=0$. For such shocks, Rankine-Hugoniot's expression \eqref{eq:RankineHugoniot} requires
\begin{align}\label{eq:HleftHright}
c=c(S_l,S_r)=\frac{h(S_l)-h(S_r)}{S_l-S_r}=0, \text{ implying } h(S_l)=h(S_r).
\end{align}
Note that this implies
\begin{align}\label{eq:SlSMSr}
S_r<S_M<S_l, \quad (\text{if } S_l>S_r).
\end{align}
Without a loss of generality we put the stationary shock at $x=0$.

\subsubsection{Case A: Connection between  imbibition and drainage  states}
Let us first consider the case
\begin{subequations}\label{eq:TwoStatesCaseI}
\begin{align}
&S(0^-,t)=S_l \text{ is in drainage state},\\
&S(0^+,t)=S_r \text{ is in imbibition state}.
\end{align}
\end{subequations}
The case of $S_l$ in imbibition and $S_r$ in drainage is symmetrical. Here, the imbibition and drainage states are inherited from the approximating solutions as before.

For \eqref{eq:TwoStatesCaseI} to be an admissible stationary shock, $\{S_l,S_r\}$ must be such that
\begin{enumerate}[label=(A\theCaseA)]
\item Relation \eqref{eq:HleftHright} is satisfied.\stepcounter{CaseA}\label{cond:Afirst}
\item Equation \eqref{eq:MainOde}, with $c(S_l,S_r)=0$, has a solution $\S(\eta)$ satisfying \eqref{eq:TWBCS}.\stepcounter{CaseA}
\item The pressure profile $\P:\R\to \R$ satisfies
\begin{align}\label{eq:PandScaseI}
\P(\eta)=\begin{cases}
\Pdr(\S(\eta)) &\text{ when } \eta<0,\\
\Pim(\S(\eta)) &\text{ when } \eta>0.\end{cases}
\end{align}\stepcounter{CaseA}\label{cond:Alast}
\end{enumerate}

Let $p_l:=\P(-\infty)$ and $p_r:=\P(+\infty)$. Then \eqref{eq:PandScaseI} implies
\begin{align}
p_l=\Pdr(S_l) \text{ and } p_r=\Pim(S_r).\label{eq:DefPlPr}
\end{align}
Putting $c(S_l,S_r)=0$ in \eqref{eq:MainOde} gives
\begin{align}\label{eq:HysMainOde}
h(\S)\P'=h(S_l)-h(S)=h(S_r)-h(S)\quad \text{ in } \R.
\end{align}
This equation implies that
\begin{align}
\P' \text{ is bounded in } \R, \text{ and thus, } \P \text{ is continuous in } \R. 
\end{align}
In particular, $\P(0^-)=\P(0^+)$,  or from \eqref{eq:PandScaseI}
\begin{align}\label{eq:PleftPright}
\Pdr(\S(0^-))=\Pim(\S(0^+)).
\end{align}
Combining \eqref{eq:PandScaseI} and \eqref{eq:HysMainOde} gives the problem 
\begin{subequations}\label{eq:ProblemEtaMin}
\begin{align}
&D^d(\S)\,\S'=h(\S)-h(S_l)\quad \text{ in } \eta<0,\\
& \S(-\infty)=S_l,
\end{align}
\end{subequations}
where $D^d(\S)=-h(\S)\frac{d \Pdr}{dS}(\S)>0$. The unique solution of \eqref{eq:ProblemEtaMin} is $\S(\eta)=S_l$ for all $\eta<0$. This follows from the fact that if $S_l>S_r$, then $S_l>S_M$ from \eqref{eq:SlSMSr} and hence $h(\S(\eta))<h(S_l)$ when $\S(\eta)>S_l$ and $h(\S(\eta))>h(S_l)$ when $S_r<\S(\eta)<S_l$. Now suppose there exists $\eta_0<0$ such that $\S(\eta_0)>S_l$. Then $h(\S(\eta_0))<h(S_l)$ implying $\S'(\eta_0)<0$. This means that $\S(\eta)>\S(\eta_0)>S_l$ for all $\eta<\eta_0$ contradicting $\S(-\infty)=S_l$. Similarly, if $S_r<\S(\eta_0)<S_l$ then $\S'(\eta_0)>0$, yielding a contradiction as well.

\begin{figure}[h!]
\begin{subfigure}{.5\textwidth}
\includegraphics[scale=.4]{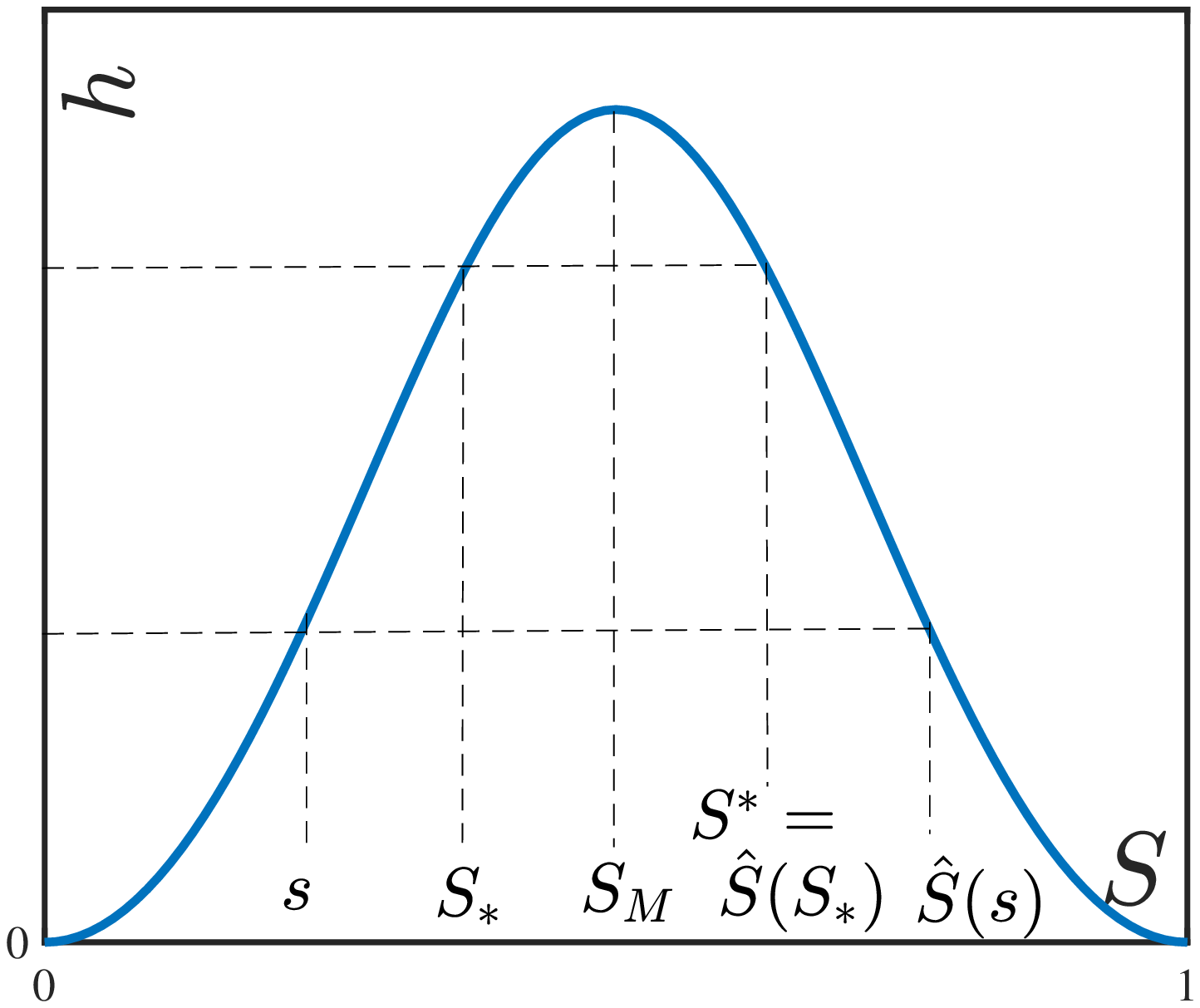}
\end{subfigure}
\begin{subfigure}{.5\textwidth}
\includegraphics[scale=.4]{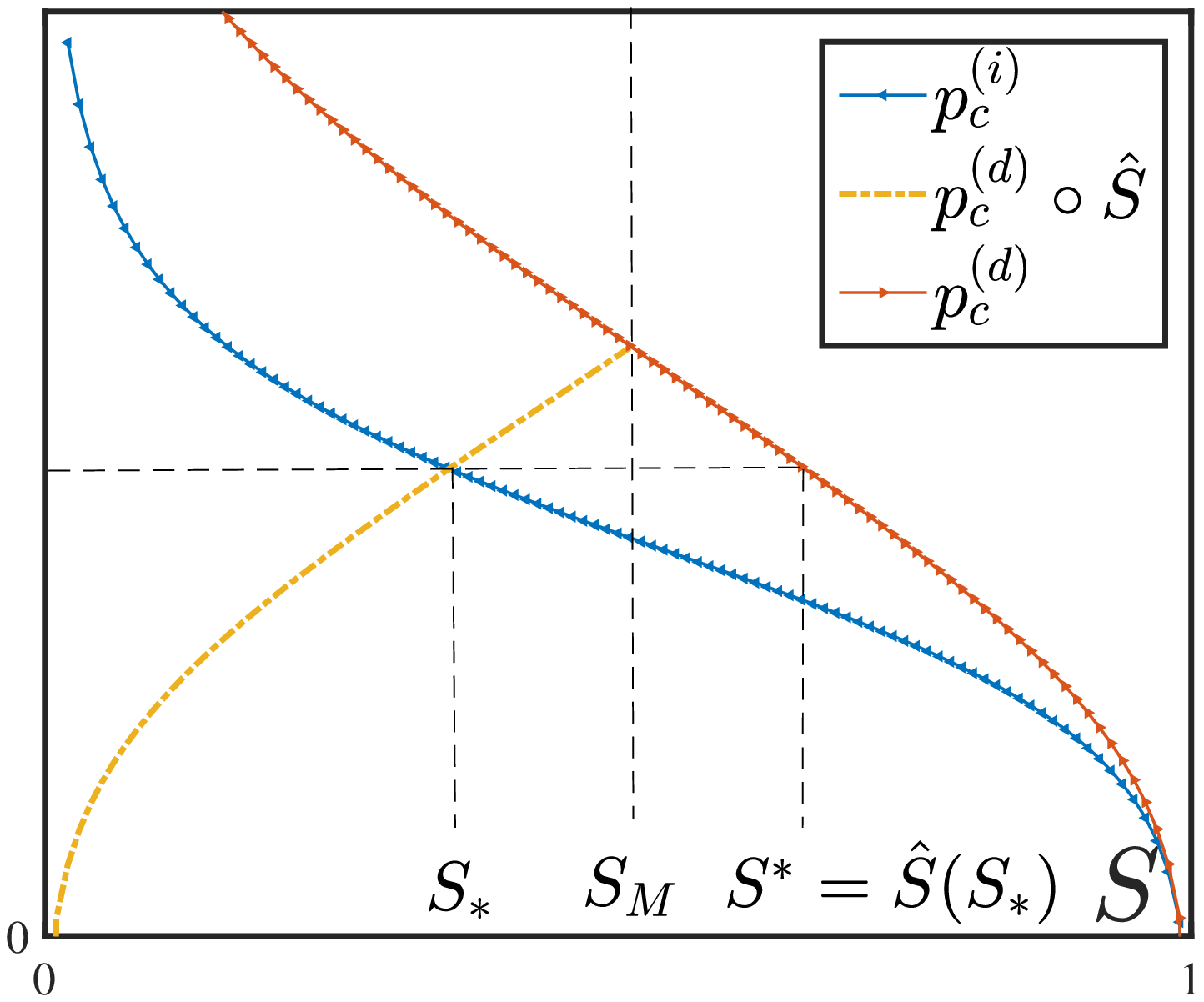}
\end{subfigure}
\caption{Construction of the pair $\{S_*,S^*\}$ defined in \eqref{eq:SupSdw}. (left) The function $\hat{S}(\cdot)$ in the $h$-$S$ diagram. (right) The composite function $\Pdr\circ \hat{S}$ and its intersection with $\Pim$.}
\label{fig:SlSr}
\end{figure}

 Repeating this reasoning for $\eta>0$, we conclude that the unique solution of \ref{cond:Afirst}--\ref{cond:Alast} is
 \begin{align}
 (\S(\eta),\P(\eta))=\begin{cases} (S_l,p_l)=(S_l,\Pdr(S_l)) &\text{ for } \eta<0,\\
(S_r,p_r)=(S_r,\Pim(S_r)) &\text{ for } \eta>0,
\end{cases}
 \end{align}
where the pair $\{S_l,S_r\}$ in addition to satisfying \eqref{eq:HleftHright}, also must satisfy $p_l=\Pdr(S_l)=\Pim(S_r)=p_r$ due to \eqref{eq:PleftPright}, i.e.,
\begin{align}\label{eq:EntropyCondI}
&\jump{h}=h(S_r)-h(S_l)=0,
&\jump{\P}:=p_r-p_l=0.
\end{align}
We show that based on the assumptions \ref{ass:Pc}--\ref{ass:h}, a unique pair $(S_l,S_r)=(S^*,S_*)$ with $0<S_*\leq S_M\leq S^*<1$ exists that satisfies \eqref{eq:EntropyCondI}. To see this, define the function $\hat{S}: [0,S_M]\to [S_M,1]$ as 
\begin{equation}
h(\hat{S}(s))=h(s).\label{eq:Def_fSr}
\end{equation}
From \ref{ass:h}, the function $\hat{S}(\cdot)$ is decreasing, continuous, $\hat{S}(0)=1$ and $\hat{S}(S_M)=S_M$, see also \Cref{fig:SlSr}  (left). Consider the strictly increasing function $\Pdr(\hat{S}(s))$ for $s\in (0,S_M)$. For $s= 0$, we have $\Pdr(\hat{S}(0))= \Pdr(1)=\Pim(1)$. Since $\Pim$ is a strictly decreasing function, for $s>0$ small enough one has $\Pdr(\hat{S}(s))< \Pim(s)$. On the other hand, $\Pdr(\hat{S}(S_M))=\Pdr(S_M)>\Pim(S_M)$. Thus, by intermediate value theorem, there exists $S_*\in (0,S_M)$ such that  $\Pdr(\hat{S}(S_*))=\Pim(S_*)$, see \Cref{fig:SlSr}  (right). Thus, we put $S^*:=\hat{S}(S_*)$.

Hence, shocks $\{S_l,S_r,0\}$ (i.e. $c(S_l,S_r)=0$), with their end states being in imbibition and drainage,  are admissible if and only if \eqref{eq:EntropyCondI} is satisfied with
\begin{subequations}\label{eq:AdmissibilityCaseI}
\begin{align}
&S_l=S^*,\; S_r=S_*,\; p_l=\Pdr(S^*),\; p_r=\Pim(S_*)\; \text{(Right state imbibition--left state drainage)},\\
& S_l=S_*,\; S_r =S^*,\; p_l=\Pim(S_*),\; p_r=\Pdr(S^*) \; \text{(Left state imbibition--right state drainage)},
\end{align}
\end{subequations}
where $\{S_*,S^*\}$ with $0<S^*\leq S_M\leq S^*<1$ is the unique solution of 
\begin{align}\label{eq:SupSdw}
 h(S_*)=h(S^*), \quad \Pim(S_*)=\Pdr(S^*).
\end{align}
Connecting an imbibition state with another imbibition state demands that \eqref{eq:EntropyCondI} be satisfied with $p_l=\Pim(S_l)$ and $p_r=\Pim(S_r)$ which has the unique solution $S_l=S_r$. The same result holds for connecting a drainage state with another drainage state. Hence, no non-trivial stationary shock ($S_l\not=S_r$) exists in these cases.
\subsubsection{Case B: One state undetermined}
\begin{figure}[h!]
\begin{subfigure}{.5\textwidth}
\includegraphics[scale=.45]{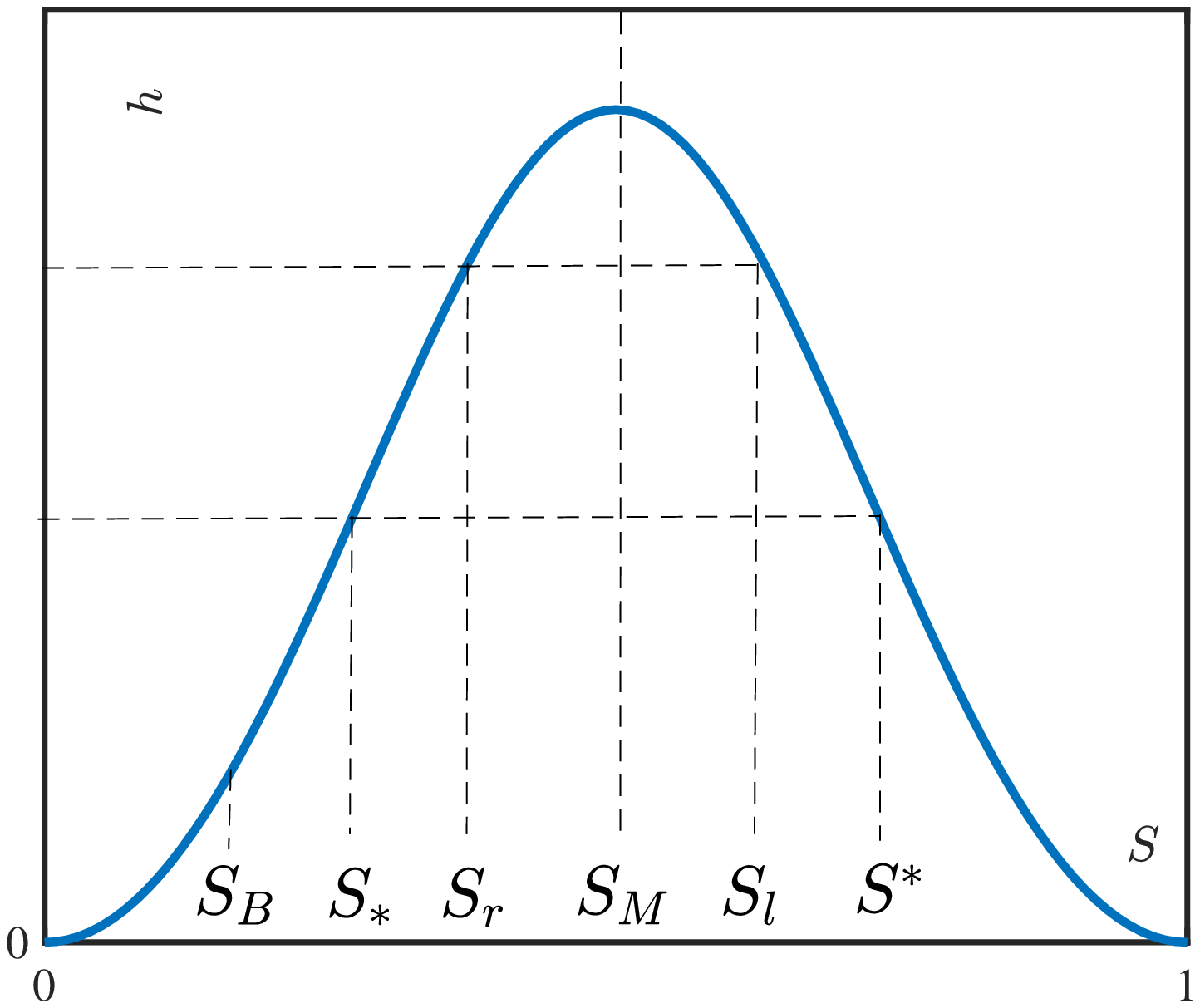}
\end{subfigure}
\begin{subfigure}{.5\textwidth}
\includegraphics[scale=.45]{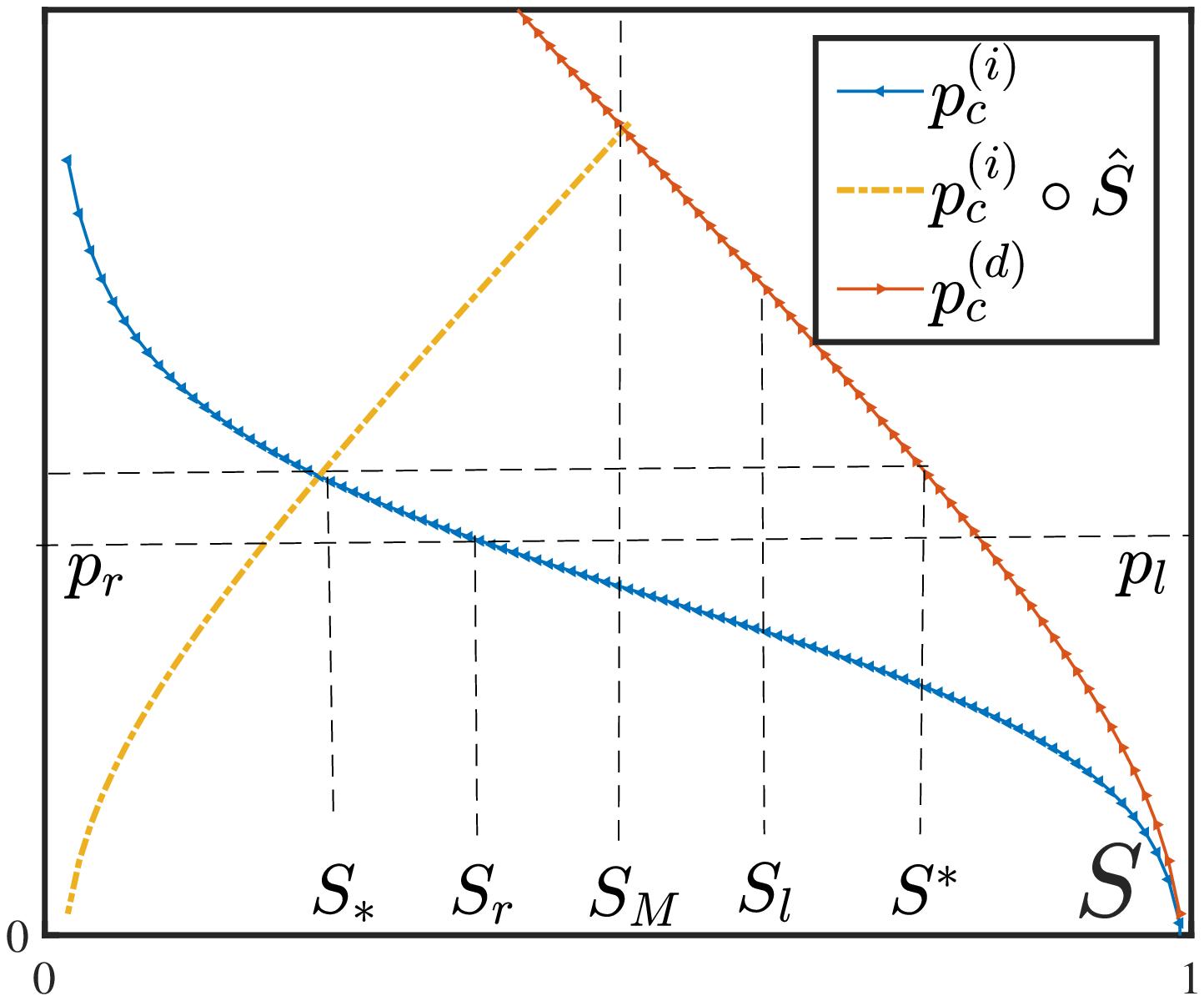}
\end{subfigure}
\caption{The saturations $S_l,\,S_r,\, S_*$ and $S_*$ in  the $h$-$S$ plane (left), and $p_c$-$S$ plane (right). The pressure $p_l=p_r=\Pim(S_r)\in [\Pim(S_l),\Pdr(S_l))$ is also shown.}\label{fig:St}
\end{figure}

It is also possible that one of the states $\{S_l,S_r\}$ is undetermined (neither in  imbibition nor in drainage). To demonstrate the admissibility for this case, we first assume that the left state is undetermined and the right state is in imbibition. Using $\p_t S= c(S_l,S_r) \S'=0$ in \eqref{eq:play-type}, we have that $\{S_l,S_r,0\}$ is an admissible stationary shock if
\begin{enumerate}[label=(B\theCaseB)]
\item Relation \eqref{eq:HleftHright} is satisfied.\stepcounter{CaseB}\label{cond:Bfirst}
\item Equation \eqref{eq:MainOde}, with $c(S_l,S_r)=0$, has a solution $\S(\eta)$ satisfying \eqref{eq:TWBCS}.\stepcounter{CaseB}
\item The pressure profile $\P:\R\to \R$ satisfies
\begin{align}\label{eq:PandScaseB}
\P(\eta)\begin{cases}
\in [\Pim(\S(\eta)),\Pdr(\S(\eta))] &\text{ when } \eta<0,\\
=\Pim(\S(\eta)) &\text{ when } \eta>0.\end{cases}
\end{align}\stepcounter{CaseA}\label{cond:Blast}
\end{enumerate}
Following the same arguments as before, we have
\begin{align}
&\P(0^-)=\P(0^+),\\
&(\S(\eta),\P(\eta))=(S_r,p_r)=(S_r,\Pim(S_r)) \text{ for all } \eta>0.
\end{align}
What remains is the $\eta<0$ problem which reads
 \begin{align}
 &\left.
\begin{matrix}
 h(\S) \P'= h(S_l)-h(\S)\\[.2em]
 \Pim(\S)\leq \P\leq \Pdr(\S)
\end{matrix}
\right\} \text{ for all } \eta<0,\\[.5em]
&\S(-\infty)=S_l,\quad \P(0^-)=p_r.
 \end{align}
Clearly, $\S=S_l$ and $\P=p_l$ is a solution provided,
\begin{subequations}\label{eq:CondOnestateUndetermined}
\begin{align}
&\jump{h}=h(S_r)-h(S_l)=0,\quad \jump{\P}= p_r-p_l=0,\\
&\Pim(S_l)< p_l=p_r =\Pim(S_r)\leq \Pdr(S_l).\label{eq:PinequalityUndetermined}
 \end{align} 
 \end{subequations}
Hence, also in this case, condition \eqref{eq:EntropyCondI} is satisfied. We show  below that additionally \eqref{eq:PinequalityUndetermined} is satisfied if
\begin{align}\label{eq:AdmissibilityCaseIIimb}
 S_r\in [S_*,S_M] \text{ and } S_l=\hat{S}(S_r)\in [S_M,S^*],
 \end{align} 
 where $\{S_*$, $S^*\}$ and the function $\hat{S}(\cdot)$ are defined in \eqref{eq:SupSdw} and \eqref{eq:Def_fSr} respectively. Indeed, from the definition of $\hat{S}(\cdot)$, if $S_r\in (S_*,S_M)$ then $h(S_l)=h(\hat{S}(S_r))=h(S_r)$ and $S_l\in (S_M,S^*)$.  From \ref{ass:Pc} it directly follows that, $\Pim(S_l)<\Pim(S_r)=p_r$. To get $\Pim(S_r)\leq \Pdr(S_l)$ we consider the strictly increasing function $p_\D(s):= \Pdr(\hat{S}(s)) -\Pim(s)$ for $s\in (0,S_M)$, see the properties of $\hat{S}(\cdot)$, \ref{ass:Pc} and \Cref{fig:SlSr,fig:St} (right). The function has a zero at $s=S_*$. This means $p_\D(S_r)\geq 0$, or $\Pdr(S_l)=\Pdr(\hat{S}(S_r))\geq \Pim(S_r)$. The construction is also shown in \Cref{fig:St}. It is straightforward to verify that this is only case when \eqref{eq:CondOnestateUndetermined} is satisfied.
 
Condition \eqref{eq:AdmissibilityCaseIIimb} serves as the admissibility criterion when the left state is undetermined and the right state is in imbibition. Considering all combinations, including undetermined-drainage, imbibition-undetermined, drainage-undetermined and undetermined-undetermined, the complete list of admissible stationary shocks $\{S_l,S_r,0\}$ with $S_l\not=S_r$ is given in \Cref{Tab:AllShocks}.

\begin{table}
\resizebox{\columnwidth}{!}{
\begin{tabular}{|l|l|l|l|l|l|}\hline \rowcolor{gray!20}
Left state &Right state &$S_l$ &$S_r$ &$p_l$ &$p_r$\\
\hline drainage &imbibition &$=S^*$ &$=S_*$ &$=\Pdr(S_l)$ &$=\Pim(S_r)$\\
\hline imbibition &drainage &$=S_*$ &$=S^*$ &$=\Pim(S_l)$ &$=\Pdr(S_r)$\\
\hline undetermined &imbibition &$\in [S_M,S^*]$ &$\in [S_*,S_M]$ &$\in (\Pim(S_l),\Pdr(S_l)]$ &$=\Pim(S_r)$\\
\hline imbibition &undetermined &$\in [S_*,S_M]$ &$\in [S_M,S^*]$ &$=\Pim(S_l)$ &$\in (\Pim(S_r),\Pdr(S_r)]$\\
\hline undetermined &drainage &$\in [S_*,S_M]$ &$\in [S_M,S^*]$ &$\in [\Pim(S_l),\Pdr(S_l))$ &$=\Pdr(S_r)$\\
\hline drainage &undetermined &$\in [S_M,S^*]$ &$\in [S_*,S_M]$ &$=\Pdr(S_l)$ &$\in [\Pim(S_r),\Pdr(S_r))$\\
\hline undetermined &undetermined &$\in [S_*,S^*]$ &$\in [S_*,S^*]$ &$\in [\Pim(S_l),\Pdr(S_l)]$ &$\in [\Pim(S_r),\Pdr(S_r)]$\\\hline
\end{tabular}}
\caption{All admissible stationary shocks $\{S_l,S_r,0\}$ with $S_l\not=S_r$. They satisfy the condition $\jump{h}=h(S_r)-h(S_l)=0,$ and $\jump{\P}= p_r-p_l=0$ where $S_l$, $S_r$, $p_l$ and $p_r$ are given in the table. The states on the left and the right of the shocks are also stated. Imbibition state  corresponds to $\Pim$ curve being  used to relate $S_l$ and $p_l$ or $S_r$ and $p_r$. Simialrly, drainage state corresponds to the use of $\Pdr$. Undetermined left state implies that $p_l\in [\Pim(S_l),\Pdr(S_l)]$ and undetermined right  state implies  $p_r\in [\Pim(S_r),\Pdr(S_r)]$.}\label{Tab:AllShocks}
\end{table}

\section{Vanishing capillarity solutions of the Buckley-Leverett equation} \label{sec:EntropySol}
We construct the solution of the Riemann problem \eqref{eq:HBL}--\eqref{eq:iniS} with the admissibility conditions from  \Cref{sec:Admissible}. A solution of \eqref{eq:HBL} is composed of constant states separated by shocks and rarefaction waves. We recall that a rarefaction wave is a smooth solution of \eqref{eq:HBL} having the form
\begin{align}
S(x,t)=r(\zeta) \text{ where } \zeta=x\slash t.
\end{align}
 Then for $r(\cdot)$ results the equation 
 \begin{align}
  \frac{d h}{d S}(r(\zeta))=\zeta \text{ for } \zeta_l<\zeta<\zeta_r,\label{eq:DefRW}
  \end{align} 
 where
 $$
r(\zeta_l)=S_l, \text{ and } r(\zeta_r)=S_r\quad  \text{ or } \quad \zeta_l=\frac{dh}{dS}(S_l), \text{ and } \zeta_r=\frac{dh}{dS}(S_r).
 $$
Rarefaction waves are well-defined if 
\begin{align}\label{eq:RWcond}
f''(s) \text{ does not change sign for all } \min\{S_l,S_r\}\leq s\leq \max\{S_l,S_r\}.
\end{align}

We show below how to construct solutions of \eqref{eq:HBL} with rarefaction waves satisfying  \eqref{eq:DefRW} along with condition \eqref{eq:RWcond}, and shocks satisfying the admissibility conditions from \Cref{sec:Admissible}. For completeness, we briefly recall the classical Buckley-Leverett construction.

\subsection{Classical construction}\label{sec:classicalEntSol}
The construction of vanishing capillarity solutions for the classical case is well established, see \cite{oleinik1957discontinuous,lefloch2002hyperbolic}. 
Recalling the definition of $\bar{S}_B$ and $\bar{S}_T$ in \eqref{eq:barS}, we identify three cases. Since the solutions for these cases will be used in \Cref{sec:HysteresisEntSol} to construct parts of the vanishing capillarity solution for pairs other than $\{S_B,S_T\}$, we also use $S_l$ and $S_r$ in the notation:
\subsubsection*{Case I: $S_T=S_l\leq S_M$:}
In this case $S_B<S_M$ since $S_B<S_T$. The vanishing capillarity solutions are
\begin{subequations}\label{eq:EntSolPlus}
\begin{align}
\begin{matrix}
\underline{\text{ If } S_l\leq \bar{S}_B}: \\\phantom{abc}
\end{matrix} \quad S(x,t)=\begin{cases}
S_l &\text{ for } x< c(S_B,S_l)\, t,\\
S_B &\text{ for } x> c(S_B,S_l)\, t,
\end{cases}
\end{align}
and with $r(\cdot)$ defined in \eqref{eq:DefRW},
\begin{align}
\begin{matrix}
\underline{\text{ If } S_l>\bar{S}_B}: \\\phantom{abc}
\end{matrix} \quad S(x,t)=\begin{cases}
S_l &\text{ for } x<h'(S_l)\, t,\\
r(x/t) &\text{ for } h'(S_l)\, t< x< c(S_B,\bar{S}_B)\, t,\\
S_B &\text{ for } x> c(S_B,\bar{S}_B) \, t.
\end{cases}
\end{align}
\end{subequations}
\subsubsection*{Case II: $S_B=S_r\geq S_M$:}
In this case $S_T>S_M$. The vanishing capillarity solutions are
\begin{subequations}\label{eq:EntSolMinus}
\begin{align}
\begin{matrix}
\underline{\text{ If } S_r\geq \bar{S}_T}: \\\phantom{abc}
\end{matrix} \quad S(x,t)=\begin{cases}
S_T &\text{ for } x< c(S_r,S_T)\, t,\\
S_r &\text{ for } x> c(S_r,S_T)\, t;
\end{cases}
\end{align}
and 
\begin{align}
\begin{matrix}
\underline{\text{ If } S_r<\bar{S}_T}: \\\phantom{abc}
\end{matrix} \quad S(x,t)=\begin{cases}
S_T &\text{ for }  x< c(\bar{S}_T,S_T)\, t,\\
r(x/t) &\text{ for } h'(S_r)\, t> x> c(\bar{S}_T,S_T)\, t,\\
S_r &\text{ for } x>h'(S_r)\, t.
\end{cases}
\end{align}
\end{subequations}
\subsubsection*{Case III: $S_B<S_M<S_T$:}
In this case we have,
\begin{equation}
S(x,t)=\begin{cases}
S_T &\text{ for } x< c(\bar{S}_T,S_T)\,t<0,\\
r(x/t) &\text{ for }  c(\bar{S}_T,S_T)\, t<x<c(S_B,\bar{S}_B)\, t,\\
S_B &\text{ for } 0<c(S_B,\bar{S}_B)\, t<x.
\end{cases} \label{eq:SolOleinik}
\end{equation}

\subsection{Construction with capillary hysteresis}\label{sec:HysteresisEntSol}

\paragraph{Vanishing capillarity solution for Case I and II:}
Observe that in Case I  of \Cref{sec:classicalEntSol} ($S_T<S_M$) the shock is in an imbibition state and the rarefaction wave part, if it exists, satisfies $\p_t S>0$. Hence, the vanishing capillarity solution as a whole is in imbibition state.
 Similarly, for Case II the entire solution is in drainage state. Consequently, based on the discussions in \Cref{sec:ShocksPlayType} and the fact that the classical solution does not depend on the specific form of $p_c$, the vanishing capillarity solutions for the capillary hysteresis model in Case I and Case II are identical to the classical ones.
 
\paragraph{Vanishing capillarity solution for Case III:}
As in the classical case, we construct a solution which is in imbibition state when $x>0$ and drainage state when $x<0$. Hence one expects to have a stationary shock at $x=0$. We restrict ourselves to the case 
\begin{align}
0<S_B<S_M<S_T<1, \quad \text{ and } \quad h(S_B)\leq h(S_T),\label{eq:RestictionSbSt}
\end{align}
the solution for $h(S_T)\leq h(S_B)$ being symmetrical. With respect to the cases A and B in \Cref{sec:ShocksPlayType} we divide the discussion in two parts.

\subsubsection{Case $S_T\geq  S^*\geq S_M$}\label{sec:Case1}
The restriction \eqref{eq:RestictionSbSt} implies that $S_B\leq S_*$.  The stationary shock at $x=0$ in this case must connect $S_*$ and $S^*$. To see this, assume the contrary. Recalling \Cref{Tab:AllShocks}, we then have
$$\lim_{x\searrow 0} S(x,t)=S_r\in (S_*,S_M] \text{ and } \lim_{x\nearrow 0} S(x,t)=\hat{S}(S_r)\in [S_M,S^*),$$
where the function $\hat{S}(\cdot)$ and $\{S_*,\,S^*\}$ are defined in \eqref{eq:Def_fSr} and  \eqref{eq:SupSdw} respectively. For $x>0$ and $t>0$, the vanishing capillarity solution of \eqref{eq:HBL} is classical, and thus, given by \eqref{eq:EntSolPlus} with $S_r\in (S_*,S_M)$. Since this is a classical solution belonging to Case I of \Cref{sec:classicalEntSol},  it is in imbibition state. Similarly, for $x<0$ and $t>0$, the solution is in drainage state. However, the only  admissible stationary shock connecting an right imbibition state to a left drainage state is $\{S^*,S_*,0\}$, see \Cref{Tab:AllShocks}. Thus, the stationary shock at $x=0$ connects $S_*$ and $S^*$.

\paragraph{Vanishing capillarity solution:} For $x>0$ and $t>0$,
\begin{subequations}\label{eq:EntSolCase1}
\begin{equation}
S(x,t) \text{ is given by } \eqref{eq:EntSolPlus} \text{ with } S_l=S_*.
\end{equation}
For $x<0$ and $t>0$,
\begin{equation}
S(x,t) \text{ is given by } \eqref{eq:EntSolMinus} \text{ with } S_r=S^*.
\end{equation}
\end{subequations}

\begin{remark}[Convergence of \eqref{eq:EntSolCase1} to the classical solutions]\label{rem:PTH2PcS}
In the absence of hysteresis, $p_c(S)=\Pim(S)=\Pdr(S)$. This means that for the standard model, $S_*=S^*=S_M$. Hence, saturation becomes continuous at $x=0$ and the vanishing capillarity solution \eqref{eq:EntSolCase1} becomes identical to \eqref{eq:SolOleinik}.
\end{remark}

\begin{figure}[h!]
\begin{subfigure}{.5\textwidth}
\begin{tikzpicture}
[xscale=3,yscale=3,domain=0:1,samples=100]
\draw[->,ultra thick] (-1,0)--(1,0) node[above, ultra thick, scale=1.5] {$x$};
\draw[->, ultra thick] (0,0)--(0,1.2) node[left,ultra thick,scale=1.5] {$t$};
  \node [below, scale=1.2] at (0,0) {0};
\draw[ultra thick, red] (0,0)--(0.7903,0.7784);
\draw[thick,dashed] (0,0)--(0.7304,0.8659);  
  \draw[thick,dashed] (0,0)--(0.6429,0.9534);
 \draw[thick,dashed] (0,0)--(0.5369,1.0117);
 \draw[thick, dashed] (0,0)--(0.4171,1.0466);
 
 \node[thick,scale=1] at (0.8641,0.8542) {$\bar{S}_B$};
\node[thick,scale=1.2] at  (0.15,0.7668) {$S_*$};
\node[thick,scale=1.2] at (0.5046,0.2187) {$S_B$};

\draw[ultra thick, red] (0,0)--(-0.9072,0.2727);
\draw[thick,dashed] (0,0)--(-0.9110,0.4266);  
  \draw[thick,dashed] (0,0)--(-0.8693,0.5338);
 \draw[thick,dashed] (0,0)--(-0.8163,0.6317);
 \draw[thick, dashed] (0,0)--(-0.7557,0.7110);
 
 \node[thick,scale=1.2] at (-0.6989,0.1) {$S_T$};
\node[thick,scale=1] at  (-0.9981,0.34) {$\bar{S}_T$};
\node[thick,scale=1.2] at (-0.3277,0.7343) {$S^*$};

\draw[ultra thick, red] (0,0)--(0,1);
\end{tikzpicture}
\end{subfigure}
\begin{subfigure}{.5\textwidth}
\begin{tikzpicture}
[xscale=3,yscale=3,domain=0:1,samples=100]
\draw[->,ultra thick] (-1,0)--(1,0) node[above, ultra thick, scale=1.5] {$x$};
\draw[->, ultra thick] (0,0)--(0,1.2) node[left,ultra thick,scale=1.5] {$t$};
  \node [below, scale=1.2] at (0,0) {0};
\draw[ultra thick, red] (0,0)--(0.8580,0.4965);

 \node[thick,scale=1.2] at (0.5663,0.1702) {$S_B$};
\node[thick,scale=1.2] at  (0.3201,0.5058) {$\check{S}_T$};
\node[thick,scale=1.2] at (-0.6610,0.4918) {$S_T$};

\draw[ultra thick, red] (0,0)--(0,1);
\end{tikzpicture}
\end{subfigure}
\caption{Vanishing capillarity solutions for the case (left) $S_T>S^*$ with the specific ordering $S_B<S_1<\bar{S}_B<S_*<S^*<\bar{S}_T<S_2<S_T$, and for the case (right) $S_M<S_T<S^*$ with the  ordering $S_B<\check{S}_T<\bar{S}_B<S_M< S_T< S^* $. The red lines represent shocks and the dashed regions are the rarefaction waves. In the left figure, there are two traveling shocks $S_B$--$\bar{S}_B$ and $S_T$--$\bar{S}_T$ and a stationary shock $S_*$--$S^*$. In the right figure, there is one shock $S_B$--$\check{S}_T$ and a stationary shock  $\check{S}_T$--$S_T$. The $S$ profiles for these two cases were shown  in \Cref{fig:StandardVsHys} and will be discussed again in \Cref{sec:Num}.}\label{fig:VanishingCapSol}
\end{figure}
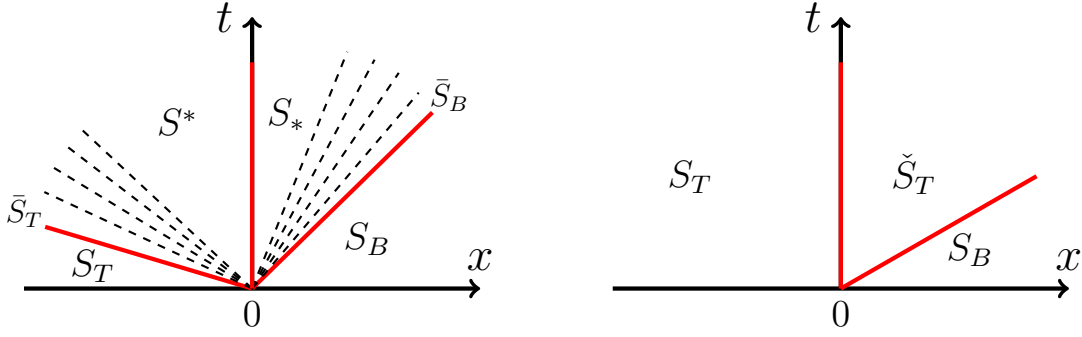

\subsubsection{Case $S^*> S_T> S_M$}\label{sec:Case2}
Let $\check{S}_T\in (S_*,S_M]$ be such that $\hat{S}(\check{S}_T)=S_T$, or in other words
\begin{align}
\check{S}_T=\hat{S}^{-1}(S_T)\in (S_*,S_M) \text{ for all } S_T\in (S_M,S^*).
\end{align}
Restriction \eqref{eq:RestictionSbSt} then implies that $S_B<\check{S}_T$.
In this case, the stationary shock $\{S^*,S_*,0\}$ is ruled out. Otherwise, since  $S_T<S^*$, the classical solution for $x<0$ and $t>0$ will be in imbibition state. However, recalling \Cref{Tab:AllShocks}, the shock $\{S^*,S_*,0\}$ is not admissible for left state in imbibition. Similarly, other possibilities are eliminated using \Cref{Tab:AllShocks}, except the following scenario:
\paragraph{Vanishing capillarity solution:} for $x>0$ and $t>0$,
\begin{subequations}\label{eq:EntSolCase2}
\begin{equation}
S(x,t) \text{ is given by } \eqref{eq:EntSolPlus} \text{ with } S_l=\check{S}_T.
\end{equation}
For $x<0$ and $t>0$,
\begin{equation}
S(x,t) =S_T \quad\text{(solution is frozen)}.
\end{equation}
\end{subequations}
The vanishing capillarity solutions for the cases $S_T>S_*$ and $S_T\in (S_M,S^*)$ are shown in \Cref{fig:VanishingCapSol}.
Observe that there is no classical counterpart of \eqref{eq:EntSolCase2}. The solution plotted in \Cref{fig:StandardVsHys} (right) for the play-type hysteresis model is of this type.

\begin{remark}[Generality of the results with respect to other hysteresis models]\label{rem:Generality}
The admissible shocks presented in  \Cref{Tab:AllShocks} and the vanishing capillarity solutions presented in \eqref{eq:EntSolCase1} and \eqref{eq:EntSolCase2} are consistent with any hysteresis model that satisfies the condition
\begin{align}\label{eq:Hysteresis}
 p\in \H, \text{ and } p=\begin{cases}
\Pim(S) &\text{ implies } \p_t S\geq 0,\\
\in [\Pim(S),\Pdr(S)] &\text{ allows } \p_t S= 0,\\
  \Pdr(S) &\text{ implies } \p_t S\leq 0.
 \end{cases}
 \end{align} 
Most of the commonly used hysteresis models, including the Lenhard-Parker model \cite{parker1987parametric} and the extended play-type model \cite{Kmitra2017} belong to this category. The results are consistent with any model satisfying  \eqref{eq:Hysteresis} since, the models only differ in the description of $(S,p)$ when $p\in (\Pim(S),\Pdr(S))$, i.e. when the hysteretic state is undetermined. As a result, if a stationary shock connects imbibition to drainage, then the set of equations \ref{cond:Afirst}--\ref{cond:Alast} are valid also for  a model satisfying \eqref{eq:Hysteresis}. Thus, the resulting shocks are unaltered. Similarly, \ref{cond:Bfirst}--\ref{cond:Blast} (in particular \ref{cond:Blast}) are consistent with describing the shock when one of the states is undetermined since in this case $\p_t S=c(S_l,S_r)\S'=0$  which is allowed by \eqref{eq:Hysteresis}. Hence, models satisfying \eqref{eq:Hysteresis} are consistent with admissible shocks listed in \Cref{Tab:AllShocks}, and consequently with the vanishing capillarity solutions derived in \Cref{sec:EntropySol}.
\end{remark}

\section{Numerical Results}\label{sec:Num}
To solve \eqref{eq:iniS}, \eqref{eq:play-type} numerically, \eqref{eq:play-type} is usually regularised in the following way:
\begin{align}\label{eq:dSdt}
 p=\begin{cases}
\Pim(S) - \t\p_t S &\text{ when } \p_t S> 0,\\
\in [\Pim(S),\Pdr(S)] &\text{ when } \p_t S= 0,\\
  \Pdr(S)- \t\p_t S &\text{ when } \p_t S< 0.
 \end{cases}
 \end{align}
The relaxation parameter $\t>0$ (dynamic capillarity coefficient) was mentioned briefly in \Cref{sec:Intro}.

We solve \eqref{eq:BL} and \eqref{eq:dSdt} in a domain $(-H,H)$ for a large $H>0$. The boundary conditions used are
\begin{align}\label{eq:NumBC}
h(S)[1+ \d\p_x p](-H,t)= h(S_T), \text{ and } p(H,t)=\Pim(S_B) \text{ for } t>0.
\end{align}
The approach for solving the system \eqref{eq:BL}, \eqref{eq:dSdt}--\eqref{eq:NumBC}, is based on rearranging \eqref{eq:dSdt} to express $\p_t S$ as a function of $S$ and $p$, and using this, to consider \eqref{eq:BL} as an elliptic equation for the pressure. Details of the numerical method are given in Section 5 of \cite{mitra2019fronts}, see also the numerical sections of \cite{mitra2018wetting,VANDUIJN2018232}. 
Cell centered finite difference method  with a uniform mesh is used for the computation with $\D x$ and $\D t$ representing the mesh and the time step sizes respectively. The following choices of parameters are made
$$
\t=0.01,\; \d=0.25,\; H= 100,\; \D x=0.01,\; \D t =0.0001.
$$
These values ensure both that our parabolic solver is converging and the capillarity solutions are sufficiently close to their hyperbolic limit. Since smaller $\d$ also implies that the profiles take longer to develop, we have optimized the value of $\d$ so that a good approximation of the developed profile is obtained in a reasonable time.

\subsection{Validation}
\begin{figure}[h!]
\begin{subfigure}{.5\textwidth}
\includegraphics[scale=.45]{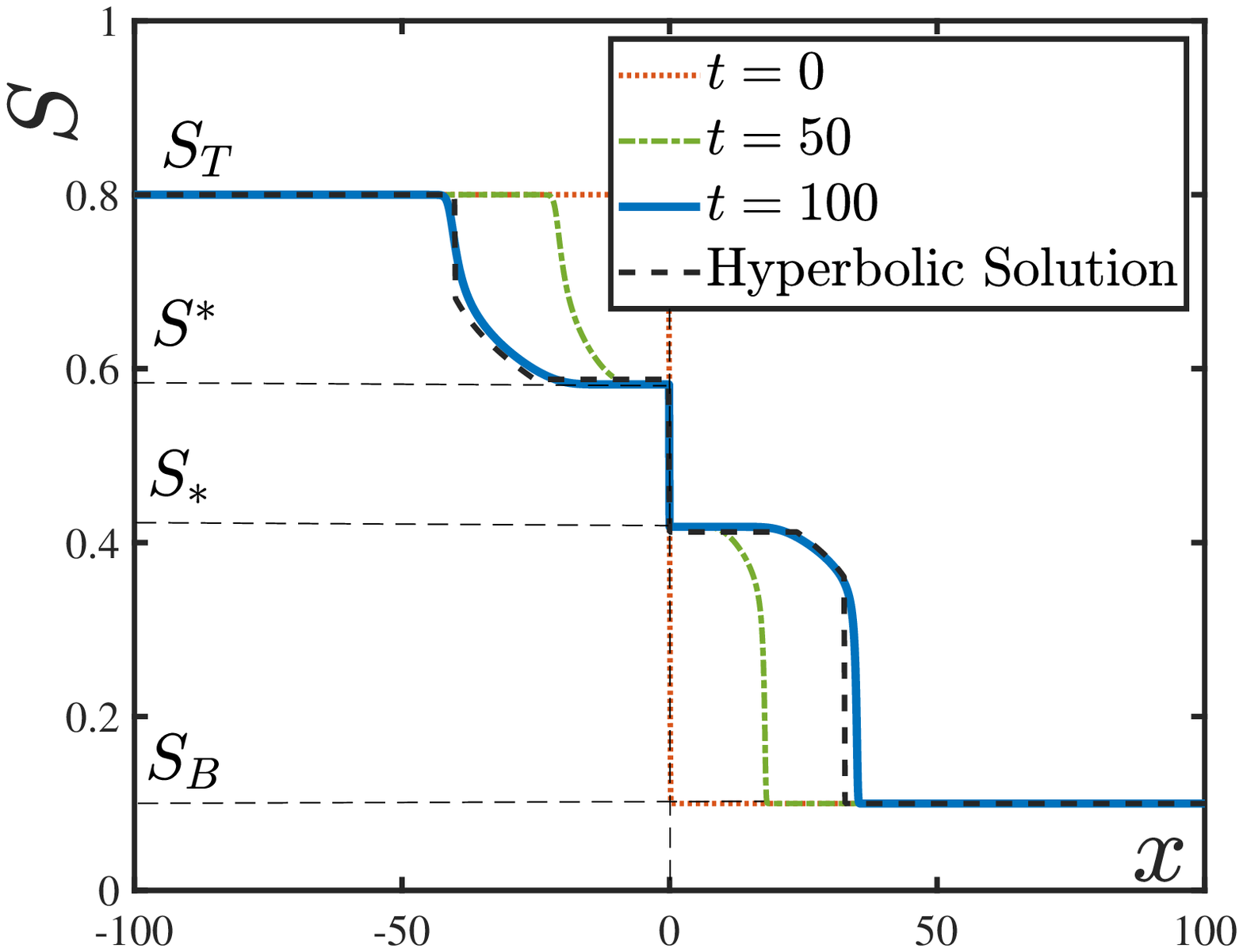}
\end{subfigure}
\begin{subfigure}{.5\textwidth}
\includegraphics[scale=.45]{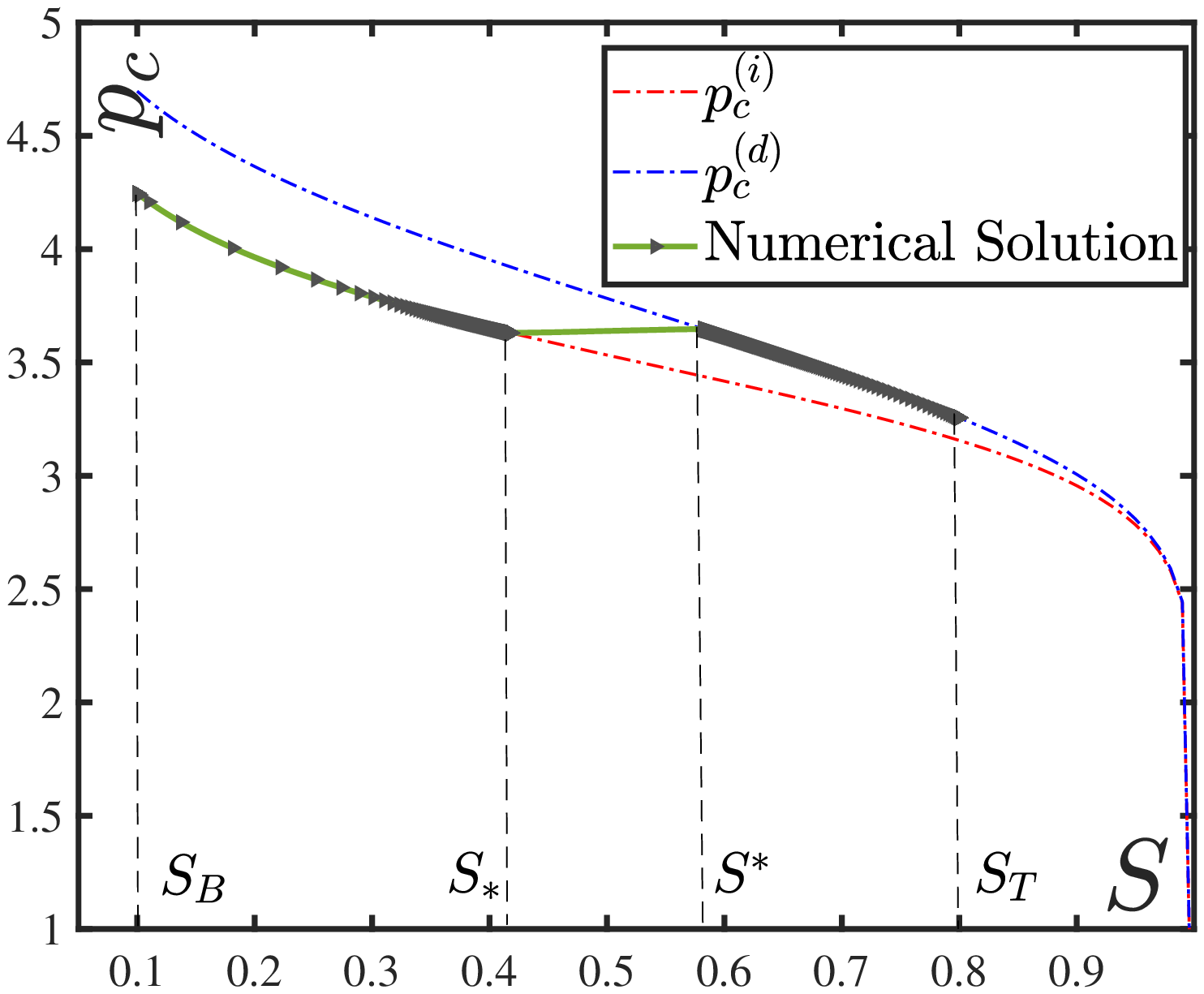}
\end{subfigure}
\caption{(left) Numerically computed capillarity solution against the vanishing capillarity solution (black dashed line) for $S_T>S^*$ given by \eqref{eq:EntSolCase1}. Here, $\d=0.25$, $S_B=0.1$ and $S_T=0.8$. The numerical solution is shown at $t=50$ and $t=100$. (right) The same numerical solution in the $p_c$-$S$ plane at $t=100$. For $x>0$, the solution is on top of $\Pim$ (imbibition state), whereas, for $x<0$ it is on top of $\Pdr$ (drainage state). A horizontal shift is observed at $x=0$ implying that the pressure is continuious across $x=0$, see \eqref{eq:EntropyCondI}.}\label{fig:Validation1}
\end{figure}
For validating our predictions of \Cref{sec:EntropySol} we use  Brooks-Corey expression \eqref{eq:hcurve} for $h$ with $M=1$. This gives $S_M=0.5$. Van Genuchten parametrization is used here to model the capillary curves. Two choices for the imbibition and drainage curves are made. At first, we take $\Pim(\cdot)$ and $\Pdr(\cdot)$ close to each other, i.e.,
$$
\Pim(S)=3.5\,(S^{-\frac{1}{q_1}}-1)^{1-q_1}  \text{ with } q_1=0.92, \text{ and } \Pdr(S)=\Pim(S) + \tfrac{1}{2}\,(1-S),
$$
see \Cref{fig:Validation1} (right). In this case, direct computation shows that $S_*=0.4121$ and $S^*=0.5879$.
We take $S_T=0.8$ and $S_B=0.1$ implying that $S_T>S^*$.
 The vanishing capillarity solution expected in this scenario is outlined in  \eqref{eq:EntSolCase1}. It consists of shocks between $S_T$--$\bar{S}_T $ and $S_B$--$\bar{S}_B$, rarefaction waves between $\bar{S}_T$--$S^*$ and $\bar{S}_B$--$S_*$ and a stationary shock from $S_*$ to $S^*$ at $x=0$, see \Cref{fig:StandardVsHys} (left) and \Cref{fig:Validation1} (left). The numerically computed solution is shown in \Cref{fig:Validation1} (left) and in the $p_c$-$S$ plane in \Cref{fig:Validation1} (right). The numerically obtained capillary  solutions are a close match to the vanishing capillarity solutions predicted. The minor differences are due to $\d>0$, $\t>0$ and numerical errors. Note that the counter-current flow for $x<0$ stems from non-monotonicity of $h$ since the flux at $x=0$ is $h(S_*)$ which is greater than the flux at the left boundary.

\begin{figure}[h!]
\begin{subfigure}{.5\textwidth}
\includegraphics[scale=.45]{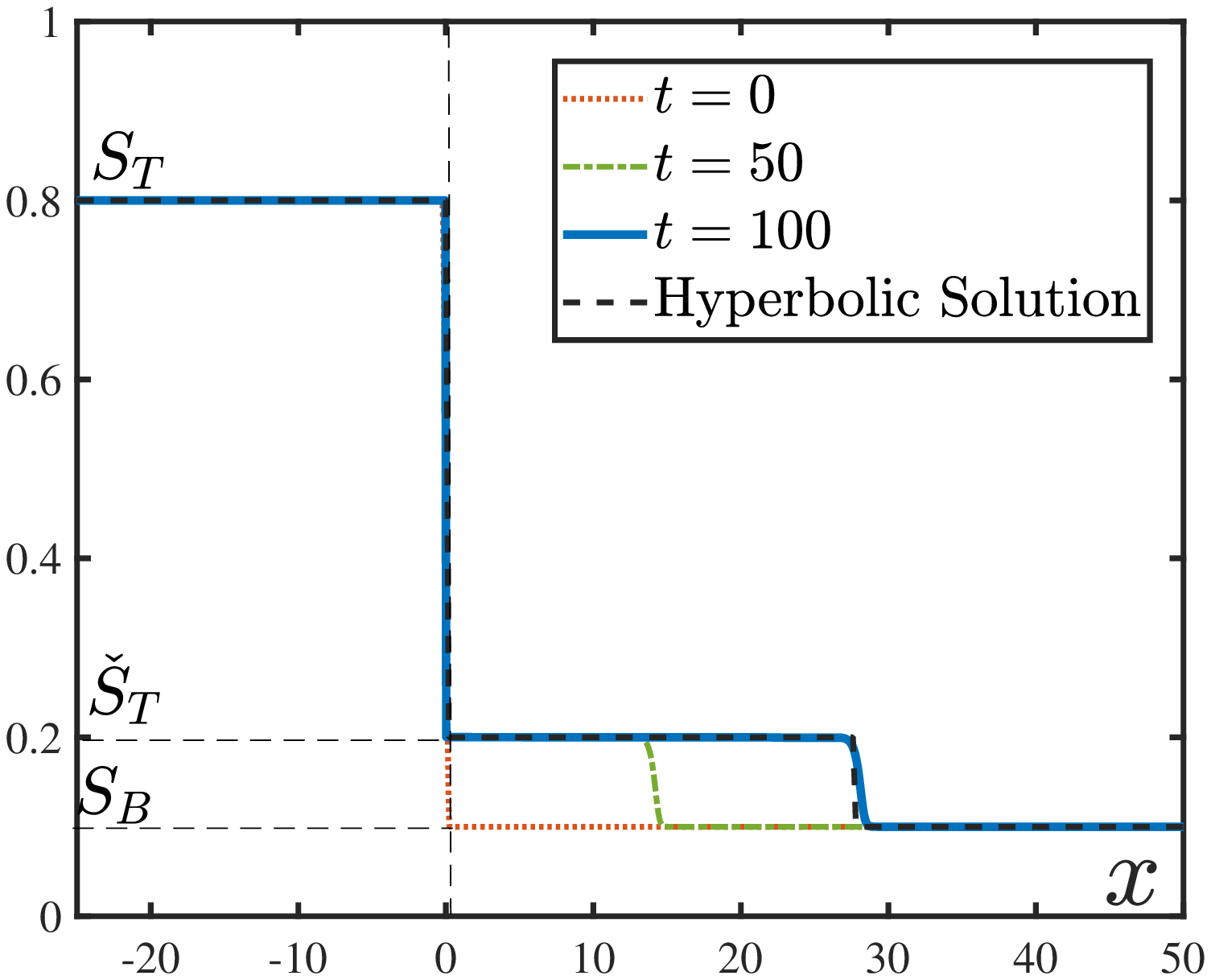}
\end{subfigure}
\begin{subfigure}{.5\textwidth}
\includegraphics[scale=.45]{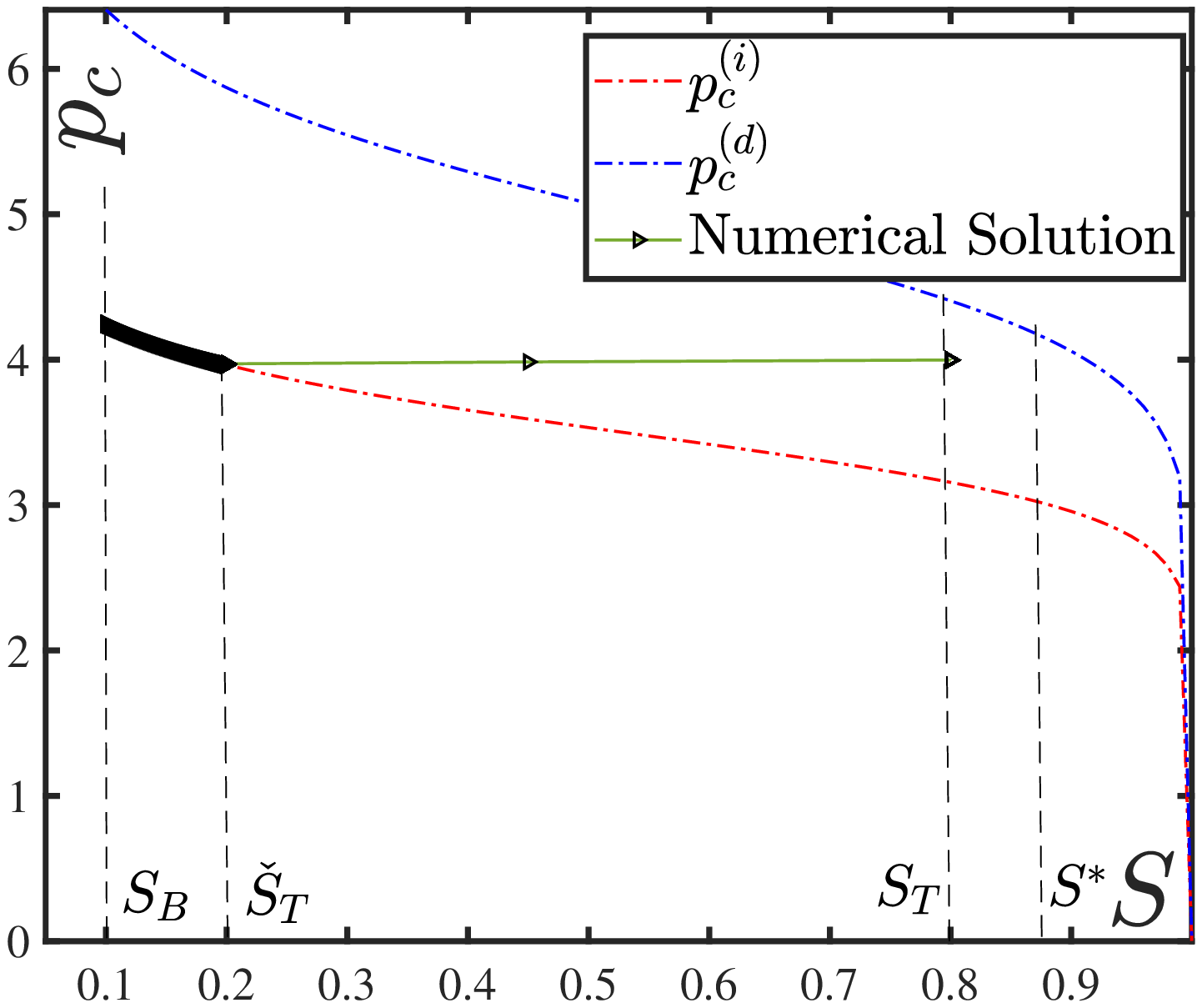}
\end{subfigure}
\caption{(left) Computed capillarity solution against the vanishing capillarity solution (black dashed line) for $S_M<S_T<S^*$ given by \eqref{eq:EntSolCase2}. Here, $\d=0.25$, $S_B=0.1$ and $S_T=0.8$. The saturation for $x<0$ does not change and $S(0,t)\approx \check{S}_T=0.2$ for all $t>0$. (right) The numerical solution in the $p_c$-$S$ plane for $t=100$. The pressure at $x=0$ is approximately $p_l=\Pim(\check{S}_T)$.}\label{fig:Validation2}
\end{figure}
Next, we take $\Pim$ same as before but $\Pdr$ as in \Cref{fig:Pch} (right):
$$
\Pdr(S)=5\, (S^{-\frac{1}{q_2}}-1)^{1-q_2}  \text{ with } q_2=0.9.
$$
These curves resemble experimentally obtained retention curves (see \cite{Morrow_Harris}) and were taken from \cite{mitra2019fronts}. With $\{S_B,S_T\}$ same as before, one has $S_T<S^*=0.8759$ in this case. Hence, recalling \Cref{sec:Case2}, particularly \eqref{eq:EntSolCase2}, the solution is frozen for $x<0$. For $x>0$, there is a shock connecting $\check{S}_T=0.2$ and $S_B$. This behaviour is mimicked by the computed capillarity solutions, see \Cref{fig:Validation2}. In this case, the total mass is still balanced since the total rate of water infiltration is $c(S_B,\check{S}_T)(\check{S}_T-S_B)= h(\check{S}_T)-h(S_B)=h(S_T)-h(S_B)$ (where $c(\cdot,\cdot)$ is the shock speed introduced in \eqref{eq:RankineHugoniot}), which is precisely equal to the difference of fluxes at the left and the right boundaries.

\section{Conclusion}
In this paper, we considered the hyperbolic Buckley-Leverett equation \eqref{eq:HBL} and its parabolic counterpart \eqref{eq:BL} for the case when the flux function $h$ is non-monotone. This occurs in gravity driven flows. In the presence of hysteresis, represented by \eqref{eq:play-type}, the hyperbolic (vanishing capillarity) limit of solutions to \eqref{eq:BL} differs from the classical solution obtained through the equilibrium expression \eqref{eq:StandardPcS}. In particular, a solution to the Riemann problem \eqref{eq:HBL}--\eqref{eq:iniS} has a stationary shock at $x=0$ if the left ($x<0$) and the right ($x>0$) states lie in the increasing and decreasing parts of the 
flux function respectively, or vice versa. The hysteretic states on the left and right become different in this case and thus, different capillary relations are used. Using travelling wave solutions, an admissibility (entropy) condition \eqref{eq:EntropyCondI} is derived for stationary shocks. The condition states that the flux function $h$ and the pressure corresponding to the hysteretic state of the system, remain continuous across the shock. It is then used to derive all admissible shocks. They are listed in \Cref{Tab:AllShocks}. The shocks are classified into two categories. The first (Case A) connects imbibition states to drainage states, and the second (Case B) has one of the states undetermined.

Depending on the values of the Riemann data with respect to characteristic points that are easily computable a-priori from $h$ and the capillary
curves, two possibilities are identified (corresponding to Case A and Case B) when the stationary shocks occur. The 
vanishing capillarity solutions for these cases are given by  \eqref{eq:EntSolCase1} and \eqref{eq:EntSolCase2}.
Interestingly, the solution \eqref{eq:EntSolCase2} remains frozen in time in one of the halves and thus,
differs significantly from the classical solution. To our knowledge, this is a novel observation. The predictions were  validated using
numerical experiments.

\section*{Acknowledgment}
K. Mitra acknowledges the support from INRIA Paris through the ERC Gatipor grant. Parts of the research for K. Mitra were also funded  by TU Dortmund University, Shell--NWO (grant 14CSER016)  and Hasselt University (grant BOF17BL04). C.J. van Duijn is supported by the Darcy Center of Eindhoven University of Technology and Utrecht University.
The authors would like to thank Prof. Sorin Pop for many fruitful discussions on the topic.


\end{document}